\documentstyle[amsthm,amsfonts,amssymb,amsgen,amsmath,amsopn,verbatim,hhline,tabularx,11pt]{article}

\theoremstyle{plain}

\theoremstyle{definition}

\newcommand{\wg}{\wedge}

\newcommand{\tZ}{{\tilde{Z}}}

\renewcommand{\P}{{\mathbb P}}

\DeclareMathOperator{\cyc}{{cyc}}

\newcommand{\eps}{{\varepsilon}}
\newcommand{\pa}{{\partial}}
\newcommand{\ga}{{\alpha}}

\newcommand{\gd}{{\delta}}

\newcommand{\CH}{{\mathcal CH}}

\begin{document}

\title{Supplement to: Goncharov's Relations in Bloch's higher Chow
Group $CH^3(F,5)$}
\author{Jianqiang Zhao\footnote{Partially supported by NSF grant DMS0139813}}
\date{}
\maketitle
\begin{center}
Department of Mathematics, Eckerd College, St. Petersburg, FL
33711
\end{center}

In this supplement we prove the admissibility of all the cycles
appearing in the paper {\em Goncharov's Relations in Bloch's
higher Chow Group $CH^3(F,5)$}. First let's recall the following
two Lemmas:

\medskip

\noindent {\bf Lemma 3.1.} (Gangl-M\"uller-Stach) {\em Let $f_i$
($i=1,2,3,5$) be rational functions and $f_4(x,y)$ be a product of
fractional linear transformations of the form
$(a_1x+b_1y+c_1)/(a_2x+b_2y+c_2)$.  We assume that all the cycles
in the lemma are {\em admissible} and write
$$Z(f_1,f_2)=[f_1,f_2,f_3,f_4,f_5]
=[f_1(x),f_2(y),f_3(x),f_4(x,y),f_5(y)]$$ if no confusion arises.

\noindent(i) If $f_4(x,y)=g(x,y)h(x,y)$ then
$$[f_1,f_2,f_3,f_4,f_5]=[f_1,f_2,f_3,g,f_5]+[f_1,f_2,f_3,h,f_5].$$

\noindent(ii) Assume that $f_1=f_2$ and that for each non-constant
solution $y=r(x)$ of $f_4(x,y)=0$ and $1/f_4(x,y)=0$ one has
$f_2(r(x))=f_2(x)$.}

\begin{quote} {\em (a) If $f_3(x)=g(x)h(x)$ then
$$[f_1,f_2,f_3,f_4,f_5]=[f_1,f_2,g,f_4,f_5]+[f_1,f_2,h,f_4,f_5].$$

(b) Similarly, if $f_5(y)=g(y)h(y)$ then
$$[f_1,f_2,f_3,f_4,f_5]=[f_1,f_2,f_3,f_4,g]+[f_1,f_2,f_3,f_4,h].$$

(c) If $f_1=f_2=gh$ and $g(r(x))=g(x)$ or $g(r(x))=h(x)$ then
\begin{equation}\label{1steq}
2 Z(f_1,f_2)=Z(g,f_2)+Z(h,f_2)+Z(f_1,g)+Z(f_1,h)
\end{equation}
and}
\begin{equation}\label{2ndeq}
Z(f_1,f_2)=Z(g,g)+Z(h,h)+Z(h,g)+Z(g,h).
\end{equation}
\end{quote}

\medskip

\noindent {\bf Lemma 3.2.} {\em Assume that $f_i$, $i=1,2,3,5$,
are rational functions of one variable and $p_4$ and $q_4$ are
rational functions of two variables. Assume that the only
non-constant solution of $p_4(x,y)=0,\infty$ is $y=x$ and the same
for $q_4(x,y)$.

(i) If $f_3=gh$ then
$$\aligned
\ [f_1,f_2,f_3, p_4, f_5] + [f_2,f_1,f_3, q_4, f_5]
=&[f_1,f_2,g, p_4, f_5]+ [f_2,f_1, g, q_4, f_5]\\
+&[f_1,f_2,h, p_4, f_5]+ [f_2,f_1, h, q_4, f_5]
\endaligned$$
if all cycles are admissible. A similar result holds if $f_5=gh$.

(ii) If $f_2=gh$ then
$$\aligned
\ [f_1,f_2,f_3, p_4, f_5] + [f_2,f_1,f_3, q_4, f_5]
=&[f_1,g, f_3,p_4, f_5]+ [g,f_1, f_3, q_4, f_5]\\
+&[f_1,h, f_3,p_4, f_5]+ [h,f_1, f_3, q_4, f_5]
\endaligned$$
if all cycles are admissible. }

\bigskip

We want to prove the following

\medskip
\noindent {\bf Theorem 4.1.} {\em Goncharov's 22 term relations
hold in $\CH^3(F,5)$: for any $a,b,c\in \P_F^1$
\begin{multline}\label{Rabc=0}
R(a,b,c)=\{-abc\}+\bigoplus_{\cyc(a,b,c)}\Bigl(
\{ca-a+1\}+\Bigl\{\frac{ca-a+1}{ca}\Bigr\}
-\Bigl\{\frac{ca-a+1}{c}\Bigr\} \\
+\Bigl\{\frac{a(bc-c+1)}{-(ca-a+1)}\Bigr\}
+\Bigl\{\frac{bc-c+1}{b(ca-a+1)}\Bigr\}
+\{c\}-\Bigl\{\frac{bc-c+1}{bc(ca-a+1)}\Bigr\}-\eta\Bigr)=0,
\end{multline}
where $\cyc(a,b,c)$ means cyclic permutations of $a,b$ and $c$,
provided that none of terms is $\{0\}$ or $\{1\}$ except for
$\eta$ (non-degeneracy condition). }

\bigskip
{\bf Step} (1). Construction of $\{k(c)\}$.

\medskip
Let $f(x)=x$, $A(x)=(ax-a+1)/a$ and $B(x)=bx-x+1$ . Let
$k(x)=B(x)/abxA(x)$ and $l(y)=1-(k(c)/k(y))$.  By definition
$$\{k(c)\}=\Bigl[x,y,1-x,1- \frac{y}{x},1- \frac{k(c)}{y}\Bigr]$$
which is easy to see as admissible. In the paper we mentioned that
for $\mu=-(ab-b+1)/a$
$$4\{k(c)\}=Z\Bigl(\frac{B}{\mu fA},\frac{B}{\mu fA}\Bigr).$$
Here for any two rational functions $f_1$ and $f_2$ of one
variable we set
$$Z(f_1,f_2)=\Bigl[f_1(x),f_2(y),1-k(x),1-\frac{k(y)}{k(x)},l(y)\Bigr].$$
We here need to prove that the following cycles
\begin{align*}
 \ &\Bigl[y,1-x,1-\frac{y}{x},1-\frac{k(c)}{y}\Bigr],\quad
 \ &\Bigl[x,1-x,1-\frac{y}{x},1-\frac{k(c)}{y}\Bigr],\quad
 \ &\Bigl[1-x,1-\frac{y}{x},1-\frac{k(c)}{y}\Bigr].
\end{align*}
are admissible and negligible.
$$\text{\boxed{\text{$Z_x:=[y,1-x,1-\frac{y}{x},1-\frac{k(c)}{y}\Bigr]$}}} $$
We have
\begin{align*}
   \pa_1^0(Z_x)&\subset \{t_3=1\},\quad
   \pa_1^\infty(Z_x)\subset \{t_4=1\},\quad
   \pa_2^\infty(Z_x)\subset \{t_3=1\},\\
   \pa_3^\infty(Z_x)&\subset \{t_2=1\},\quad
   \pa_4^\infty(Z_x)\subset \{t_3=1\},
\end{align*}
and
$$\pa_2^0(Z_x)=\pa_3^0(Z_x)=\Bigl[x,1-x,1-\frac{k(c)}{x}\Bigr],\quad
   \pa_4^0(Z_x)=\Bigl[k(c),1-x,1-\frac{k(c)}{x}\Bigr],$$
which are both admissible because
\begin{equation}\label{kcne1}
1-k(c)=\frac{(c-1)(1+abc)}{abcA(c)}\ne 0.
\end{equation}

$$\text{\boxed{\text{$Z_y:=[x,1-x,1-\frac{y}{x},1-\frac{k(c)}{y}\Bigr]$}}} $$
We have
\begin{align*}
   \pa_1^0(Z_y)&\subset \{t_2=1\},\quad
   \pa_1^\infty(Z_y)\subset \{t_3=1\},\quad
   \pa_2^0(Z_y)\subset \{t_1=1\},\\
   \pa_2^\infty(Z_y)&\subset \{t_3=1\},\quad
   \pa_3^\infty(Z_y)\subset \{t_2=1\},\quad
   \pa_4^\infty(Z_y)\subset \{t_3=1\},
\end{align*}
and
$$\pa_3^0(Z_y)=\pa_4^0(Z_y)=\Bigl[x,1-x,1-\frac{k(c)}{x}\Bigr]=\pa_2^0(Z_x)$$
which is admissible.

$$\text{\boxed{\text{$Z_{x,y}:=[1-x,1-\frac{y}{x},1-\frac{k(c)}{y}\Bigr]$}}} $$
We have
\begin{align*}
   \pa_1^\infty(Z_y)\subset \{t_2=1\},\quad
   \pa_2^\infty(Z_y)&\subset \{t_1=1\},\quad
   \pa_3^\infty(Z_y)\subset \{t_2=1\},
\end{align*}
and
$$
\pa_1^0(Z_x)=\pa_2^0(Z_y)=\pa_3^0(Z_y)=\Bigl[1-y,1-\frac{k(c)}{y}\Bigr]$$
which is admissible by \eqref{kcne1}.

\medskip
\noindent {\bf Step} (2).  The key reparametrization and a simple
expression of $\{k(c)\}$.

\medskip
Applying Lemma 3.1(ii) we see that
\begin{multline}\label{AA}
4\{k(c)\} = Z\Bigl(\frac{\mu fA}{B}, \frac{\mu fA}{B}\Bigr)
    =Z(A,A)+Z\Bigl(\frac{\mu f}{B},A\Bigr)
    +Z\Bigl(A,\frac{\mu f}{B}\Bigr)
    +Z\Bigl(\frac {\mu f}{B},\frac{\mu f}{B}\Bigr) \\
= Z(A,A)+\rho_x Z(A,A)+\rho_y Z(A,A)+\rho_{x,y}Z(A,A)=4Z(A,A).
\end{multline}
 We only need to show that the following cycle is admissible:
$$\text{\boxed{\text{$Z_A:=Z(A,A)=\Bigl[A(x),A(y),1-k(x),1-\frac{k(y)}{k(x)}, l(y)\Bigr]$}}} $$
Note that
\begin{align} \label{1-kx}
1-k(x)=&\frac{(x-1)(1+abx)}{abxA(x)},\\
1-\frac{k(y)}{k(x)}=&\frac{(y-x)(yB(x)+A(x))}{yA(y)B(x)}=\frac{(y-x)(xB(y)+A(y))}{yA(y)B(x)}.
\label{1-kyx}
\end{align}
We have
\begin{align*}
    \pa_1^0(Z_A)&\subset \{t_4=1\},\quad
   \pa_1^\infty(Z_A)\subset \{t_3=1\},
   \pa_2^0(Z_A)\subset \{t_5=1\},\quad
   \pa_2^\infty(Z_A)\subset \{t_4=1\},\\
   \pa_3^\infty(Z_A)&\subset \{t_4=1\},\quad
   \pa_4^\infty(Z_A)\subset \{t_3=1\}\cup \{t_5=1\},\quad
   \pa_5^\infty(Z_A)\subset \{t_4=1\},\\
    \pa_3^0(Z_A)&=\Bigl[\frac{1}{a},A(y),1-k(y),l(y)\Bigr]
    +\Bigl[A\Big(\frac{-1}{ab}\Big),A(y),1-k(y),l(y)\Bigr],\\
   \pa_4^0(Z_A)&=\Bigl[A(y),A(y),1-k(y),l(y)\Bigr]
   +\Bigl[\frac{\mu y}{B(y)},A(y),1-k(y),l(y)\Bigr],\\
   \pa_5^0(Z_A)&=\Bigl[A(x),A(c),1-k(x),l(x)\Bigr]
   +\Bigl[A(x),A(y_2),1-k(x),l(x)\Bigr],
\end{align*}
where the last equation comes from the two solutions of $l(y)=0$:
\begin{equation}\label{equ:ly=0}
y_1=c \quad\text{ and }  \quad
y_2=-\frac{ac-a+1}{a(bc-c+1)}=-\frac{A(c)}{B(c)}=\rho_c(c).
\end{equation}
By non-degeneracy assumption and
\begin{align}\label{Ay2}
A(y_2)= &\rho_c(A(c))=c\mu /B(c)\ne 0,\infty, \\
B(y_2)= &\rho_c(B(c))=-\mu/B(c) \ne 0,\infty, \label{By2}
\end{align}
it suffices to show the following cycles are admissible:
\begin{alignat*}{3}
    L&:=\Bigl[A(y),1-k(y),l(y)\Bigr],\quad&
    L'&:=\Bigl[A(y),A(y),1-k(y),l(y)\Bigr],\\
    L''&:=\Bigl[\frac{\mu y}{B(y)},A(y),1-k(y),l(y)\Bigr].
\end{alignat*}

\begin{itemize}
\item $L$ is admissible. Because $l(y)=1-yB(c)A(y)/cA(c)B(y)$ we
have
$$\pa_1^0(L)\subset \{t_3=1\},\quad
 \pa_1^\infty(L)\subset\{t_2=1\},\quad
 \pa_2^\infty(L)\subset\{t_3=1\},\quad
 \pa_3^\infty(L)\subset\{t_2=1\}.
$$
Moreover, by non-degeneracy assumption we see that (note that
$k(1)=k(-1/ab)=1$ by \eqref{1-kx})
\begin{align}
A(1)=&\frac{1}{a}\ne 0,\infty,\quad
l(1)=1-k(c)=\frac{(c-1)(1+abc)}{bc(ca-a+1)}\ne
0,\infty,  \label{l1}\\
A\Big(\frac{-1}{ab}\Big)=&\frac{\mu}{b}\ne 0,\infty,\quad
l\Big(\frac{-1}{ab}\Big)= 1-k(c)\ne 0,\infty, \label{lab}\\
aby_2+1=&\frac{(1-c)(ab-b+1)}{bc-c+1} \ne 0,\infty. \label{aby2}
\end{align}
Thus both $\pa_2^0(L)=[A(1),l(1)]+[A(-1/ab),l(-1/ab)]$ and
$\pa_3^0(L)=[A(c),1-k(c)]+[A(y_2),1-k(c)]$ are clearly admissible
by non-degeneracy assumption, \eqref{l1}, \eqref{lab}, and
\eqref{aby2}.

\item $L'$ is admissible. This follows from the above proof for
$L$.

\item $L''$ is admissible. This also follows from the proof for
$L$ because $\mu y/B(y)\ne 0,\infty$ when $y=1,-1/ab,c, y_2$ by
\eqref{By2}.

\end{itemize}

\medskip
\noindent {\bf Step} (3). Some admissible cycles for decomposition
of $\{k(c)\}$.

\medskip

Define the following cycles
\begin{align*}
Z_1(A,A)=&\Bigl[\frac{(b-1)A(x)}{\mu},\frac{(b-1)A(y)}{\mu},
        \frac{x-1}{x}, \frac{y-x}{A(y)}, l(y)\Bigr]\\
Z_1=&\Bigl[A(x),A(y),\frac{x-1}{x}, \frac{y-x}{A(y)}, l(y)\Bigr],\\
Z_2(A,A)=&\Bigl[\frac{(b-1)A(x)}{\mu},\frac{(b-1)A(y)}{\mu},
        \frac{x-1}{x}, \Bigl(\frac{A(y)}{y}\Bigr)
    \Bigl(1-\frac{\mu x}{A(y)B(x)}\Bigr),l(y)\Bigr]\\
Z_2=&\Bigl[\frac{(b-1)A(x)}{\mu},\frac{(b-1)A(y)}{\mu},
        \frac{x-1}{x}, 1-\frac{\mu x}{A(y)B(x)},l(y)
\Bigr],\\
Z_3(A,A) =&\Bigl[\frac{(b-1)A(x)}{\mu},\frac{(b-1)A(y)}{\mu},
        \frac{abx+1}{abA(x)}, \frac{y-x}{A(y)},l(y)\Bigr]\\
Z_3=&\Bigl[A(x),A(y),\frac{abx+1}{abA(x)},
        \frac{y-x}{A(y)}, l(y)\Bigr],\\
Z_4(A,A)=&\Bigl[\frac{(b-1)A(x)}{\mu},\frac{(b-1)A(y)}{\mu},
        \frac{abx+1}{abA(x)}, \Bigl(\frac{A(y)}{y}\Bigr)
        \Bigl(1-\frac{\mu x}{A(y)B(x)}\Bigr),l(y)\Bigr]\\
Z_4=&\Bigl[\frac{(b-1)A(x)}{\mu},\frac{(b-1)A(y)}{\mu},
        \frac{abx+1}{aA(x)}, \Bigl(\frac{A(y)}{y}\Bigr)
        \Bigl(1-\frac{\mu x}{A(y)B(x)}\Bigr),l(y)\Bigr].
\end{align*}

\noindent{\bf Claim 1}. Modulo admissible and negligible cycles
the two cycles $Z_1(A,A)$ and $Z_1$ are the same and both
admissible.

$$\text{\boxed{\text{$Z_1=\Bigl[A(x),A(y),\frac{x-1}{x}, \frac{y-x}{A(y)}, l(y)\Bigr]$}}} $$

On $\pa_1^0(Z_1)$ we have $x=1-a^{-1}$ and $(y-x)/A(y)=1$. By
similar argument we get
\begin{align*}
 \pa_1^0(Z_1)\subset &\{t_4=1\},\quad\pa_1^\infty(Z_1)\subset\{t_3=1\},\quad
 \pa_2^0(Z_1)\subset \{t_5=1\}, \\
 \pa_2^\infty(Z_1) \subset & \{t_4=1\},\quad
  \pa_4^\infty(Z_1)\subset\{t_5=1\}.
\end{align*}
So we still need to show the following cycles are admissible:
\begin{align*}
\pa_3^0(Z_1)=&\Big[\frac{1}{a}, A(y),\frac{y-1}{A(y)},
l(y)\Big],\qquad
\pa_3^\infty(Z_1)=\Big[\frac{1-a}{a},A(y),\frac{y}{A(y)}, l(y)\Big],\\
\pa_4^0(Z_1)=&\Big[A(y),A(y),\frac{y-1}{y}, l(y)\Big],\quad
\pa_5^\infty(Z_1)= \Bigl[A(x),\frac{\mu}{b-1},
 \frac{x-1}{x},\frac{B(x)}{-\mu} \Bigr],\\
\pa_5^0(Z_1)=&\Bigl[A(x),A(c),\frac{x-1}{x},\frac{c-x}{A(c)}\Bigr]
+\Bigl[A(x),A(y_2),\frac{x-1}{x},\frac{y_2-x}{A(y_2)}\Bigr].
\end{align*}
So we still need to show the following cycles are admissible:
\begin{alignat*}{3}
T:=&[A(y),(y-1)/A(y), l(y)],\quad &
U:=&[A(y),y/A(y), l(y)],\\
V:=&[A(y),A(y),1-y^{-1}, l(y)],\quad &
W:=&\Bigl[A(x),\frac{x-1}{x},\frac{B(x)}{-\mu} \Bigr],\\
X:=&\Bigl[A(x),\frac{x-1}{x},\frac{c-x}{A(c)} \Bigr],\quad&
Y:=&\Bigl[A(x),\frac{x-1}{x}, \frac{y_2-x}{A(y_2)}\Bigr].
\end{alignat*}

\begin{itemize}
\item $T$ is admissible. Because
$l(y)=1-\frac{yB(c)A(y)}{cA(c)B(y)}$ we have
$$\pa_1^0(T)\subset \{t_3=1\},\quad
 \pa_1^\infty(T)\subset\{t_2=1\},\quad
 \pa_2^\infty(T)\subset\{t_3=1\}.
$$
Moreover $\pa_2^0(T)=[1/a,l(1)]$ and
$\pa_3^\infty(T)=[\mu/(b-1),-b/\mu]$ are clearly admissible by
non-degeneracy assumption and \eqref{l1}. Lastly, from the two
solutions of $l(y)=0$ in \eqref{equ:ly=0} and \eqref{Ay2} we get
$$\pa_3^0(T)=\Bigl[A(c), \frac{c-1}{A(c)}\Bigr]+
 \Bigl[A(y_2), \frac{y_2-1}{A(y_2)} \Bigr]$$
which is admissible by non-degeneracy assumption.

\item $U$ is admissible. Similar to $T$ we have
$$\pa_1^0(U)\subset \{t_3=1\},\quad
 \pa_1^\infty(U)\subset\{t_2=1\},\quad
 \pa_2^0(U)\subset \{t_3=1\}, \quad
 \pa_2^\infty(U)\subset\{t_3=1\}.
$$
Moreover, $\pa_3^\infty(U)=[\mu/(b-1),a/(ab-b+1)]$ is clearly
admissible by non-degeneracy assumption. Lastly, from
\eqref{equ:ly=0} and \eqref{Ay2} we get
$$\pa_3^0(U)=\Bigl[A(c), \frac{c}{A(c)}\Bigr]+
\Bigl[A(y_2), \frac{y_2}{A(y_2)} \Bigr]$$ which is admissible by
non-degeneracy assumption.

\item $V$ is admissible.  First it's easy to see that
\begin{align*}
 \pa_1^0(V) &\subset\{t_4=1\},\quad
 \pa_1^\infty(V)\subset\{t_3=1\},\quad
 \pa_2^0(V)  \subset \{t_4=1\}, \quad
 \pa_2^\infty(V) \subset  \{t_3=1\},\\
 \pa_3^0(V)&=[1/a,1/a,l(1)],\quad
 \pa_3^\infty(V)\subset\{t_4=1\},\quad
\pa_4^\infty(V)=\Bigl[\frac{\mu}{b-1},\frac{\mu}{b-1},b\Bigr].
\end{align*}
 From \eqref{l1} both cycles are clearly
admissible by non-degeneracy assumption. Lastly, from the two
solutions of $l(y)=0$ in \eqref{equ:ly=0} we get
$$\pa_4^0(V)=\Bigl[A(c),A(c),\frac{c-1}{c}\Bigr]+
 \Bigl[A(y_2),A(y_2),\frac{y_2-1}{y_2}\Bigr]$$
which is admissible by \eqref{Ay2} and
\begin{equation}\label{y2nondeg}
\frac{y_2-1}{y_2}=\frac{1+abc}{ac-a+1}\ne \infty,0.
\end{equation}

\item $W$ is admissible. We can compute as follows:
\begin{alignat*}{3}
\pa_1^0(W)&\subset\{t_3=1\},\quad&
\pa_1^\infty(W)\subset\{t_2=1\},\quad&
\pa_2^0(W)=\Big[\frac{1}{a},\frac{-b}{\mu}\Big], \\
\pa_2^\infty(W)&=\Big[\frac{1-a}{a},\frac{-1}{\mu}\Big],\quad&
 \pa_3^0(W)=\Big[-\mu,b\Big],\quad & \pa_3^\infty(W)\subset\{t_2=1\}.
\end{alignat*}
All the cycles above are clearly admissible.

\item $X$ is admissible. Similar to $W$ we have
\begin{alignat*}{5}
 \pa_1^0(X)&\subset \{t_3=1\},\quad&
 \pa_1^\infty(X)&\subset\{t_2=1\},\quad&
 \pa_2^0(X)&=\Bigl[\frac{1}{a},\frac{c-1}{A(c)}\Bigr],\\
 \pa_2^\infty(X)&=\Bigl[\frac{1-a}{a},\frac{c}{A(c)}\Bigr], \quad&
 \pa_3^0(X)&=\Bigl[A(c), \frac{c-1}{c}\Bigr],\quad&
\pa_3^\infty(X)&\subset\{t_2=1\}.
\end{alignat*}
All the cycles above are clearly admissible.

\item $Y$ is admissible.  Similar to the above we get
\begin{alignat*}{5}
\pa_1^0(Y) & \subset \{t_3=1\},\quad &
 \pa_1^\infty(Y)& \subset\{t_2=1\},\quad &
 \pa_2^0(Y)& =\Bigl[\frac{1}{a},\frac{y_2-1}{A(y_2)}\Bigr], \\
 \pa_2^\infty(Y)& =\Bigl[\frac{1-a}{a},\frac{y_2}{A(y_2)}\Bigr],\quad &
 \pa_3^0(Y)& =\Bigl[A(y_2), \frac{y_2-1}{y_2}\Bigr], \quad & \pa_3^\infty(Y) & \subset\{t_2=1\}.
\end{alignat*}
By \eqref{equ:ly=0}, \eqref{y2nondeg} and \eqref{Ay2} all the
cycles above are clearly admissible. This concludes the proof that
$Z_1$ is an admissible cycle.
\end{itemize}

$$\text{\boxed{\text{$ Z_1(A,A)=\Bigl[\frac{(b-1)A(x)}{\mu},
 \frac{(b-1)A(y)}{\mu},\frac{x-1}{x}, \frac{y-x}{A(y)}, l(y)\Bigr]$}}} $$
Throughout the proof that $Z_1$ is an admissible cycle we never
use the hyperplane $\{t_1=1\}$ and moreover $(b-1)/\mu\ne
0,\infty$ by non-degeneracy assumption. Therefore we can use
exactly the same proof to show the admissibility of $Z_1(A,A)$.

$$\text{\boxed{\text{$Z_{11}=\Bigl[\frac{b-1}{\mu},A(y),\frac{x-1}{x},
 \frac{y-x}{A(y)}, l(y)\Bigr]\in C^1(F,1)\wg C^2(F,4)$}}} $$
Let
$$Z'_{11}=\Bigl[A(y),\frac{x-1}{x}, \frac{y-x}{A(y)}, l(y)\Bigr]$$

It is easy to see that
\begin{align*}
 \pa_1^0(Z'_{11})&\subset \{t_4=1\},\qquad\qquad\qquad\quad\ \pa_1^\infty(Z'_{11})\subset\{t_3=1\},\\
 \pa_2^0(Z'_{11})&=\Bigl[A(y), \frac{y-1}{A(y)}, l(y)\Bigr]=T,\quad \
 \pa_2^\infty(Z'_{11})=\Bigl[A(y), \frac{y}{A(y)}, l(y)\Bigr]=U,\\
 \pa_3^0(Z'_{11})&=\Bigl[A(y), \frac{y-1}{y}, l(y)\Bigr]=:V', \quad
\pa_3^\infty(Z'_{11})\subset \{t_4=1\}, \\
\pa_4^0(Z'_{11})&=\Bigl[A(c),\frac{x-1}{x},
\frac{c-x}{A(c)}\Bigr]+
\Bigl[A(y_2),\frac{x-1}{x}, \frac{y_2-x}{A(y_2)}\Bigr]=:X'+Y',\\
\pa_4^\infty(Z'_{11})&=\Bigl[\frac{\mu}{b-1}, \frac{x-1}{x},
\frac{B(x)}{-\mu}\Bigr]=:W'.
\end{align*}
The admissibility of $V',X',Y'$ and $W'$ follows from the proof of
that of $V,X,Y$ and $W$, respectively.

$$\text{\boxed{\text{$Z_{12}=\Bigl[A(x),\frac{b-1}{\mu},\frac{x-1}{x},
 \frac{y-x}{A(y)}, l(y)\Bigr]\in C^1(F,1)\wg C^2(F,4)$}}} $$
Let
$$Z'_{12}=\Bigl[A(x),\frac{x-1}{x}, \frac{y-x}{A(y)}, l(y)\Bigr] $$
It is easy to see that
\begin{align*}
 \pa_1^0(Z'_{12})&\subset \{t_3=1\},\quad \pa_1^\infty(Z'_{12})\subset\{t_2=1\},\\
 \pa_2^0(Z'_{12})&=\Bigl[\frac{1}{a}, \frac{y-1}{A(y)},l(y)\Bigr]=T',\quad
 \pa_2^\infty(Z'_{12})=\Bigl[\frac{1-a}{a}, \frac{y}{A(y)}, l(y)\Bigr]=U',\\
 \pa_3^0(Z'_{12})&=\Bigl[A(y), \frac{y-1}{y}, l(y)\Bigr]=V', \quad
\pa_3^\infty(Z'_{12})\subset \{t_4=1\}, \\
\pa_4^0(Z'_{12})&=\Bigl[A(x),\frac{x-1}{x},
\frac{c-x}{A(c)}\Bigr]+
\Bigl[A(x),\frac{x-1}{x}, \frac{y_2-x}{A(y_2)}\Bigr]=X+Y,\\
\pa_4^\infty(Z'_{12})&=\Bigl[A(x), \frac{x-1}{x},
\frac{B(x)}{-\mu}\Bigr]=W.
\end{align*}
The admissibility of $T',U'$ and $V'$ follows from the proof of
that of $T,U$ and $V$, respectively.

$$\text{\boxed{\text{$Z_{13}=\Bigl[\frac{b-1}{\mu},\frac{b-1}{\mu},
 \frac{x-1}{x}, \frac{y-x}{A(y)}, l(y)\Bigr]\in C^1(F,1)\wg C^1(F,1)\wg C^1(F,3)$}}} $$
Let
$$Z'_{13}=\Bigl[\frac{x-1}{x}, \frac{y-x}{A(y)}, l(y)\Bigr]$$
It is easy to see that $\pa_2^\infty(Z'_{13})\subset \{t_3=1\}$
and
\begin{align*}
 \pa_1^0(Z'_{13})&=\Bigl[\frac{y-1}{A(y)}, l(y)\Bigr]=:T'',\quad \
 \pa_1^\infty(Z'_{13})=\Bigl[\frac{y}{A(y)}, l(y)\Bigr]=:U'',\\
 \pa_2^0(Z'_{13})&=\Bigl[\frac{y-1}{y}, l(y)\Bigr]=:V'',\quad
\pa_3^\infty(Z'_{13})=\Bigl[\frac{x-1}{x},
\frac{B(x)}{-\mu}\Bigr]=:W'',\\
\pa_3^0(Z'_{13})&=\Bigl[\frac{x-1}{x}, \frac{c-x}{A(c)}\Bigr]+
\Bigl[\frac{x-1}{x}, \frac{y_2-x}{A(y_2)}\Bigr]=:X''+Y''.
\end{align*}
The admissibility of $T'',U'',V'',X'',Y''$ and $W''$ is easy to
check. It also follows from the proof of the admissibility of
$T,U,V,X,Y$ and $W$, respectively.

All the above justifies the use of Lemma 3.1(ii)(c)(2) to get:
$$Z_1(A,A)=Z_1+Z_{11}+Z_{12}+Z_{13}.$$

\medskip

\noindent{\bf Claim 2}. Modulo admissible and negligible cycles
the two cycles $Z_2(A,A)$ and $Z_2$ are the same and both
admissible.

$$\text{\boxed{\text{$Z_2=\Bigl[\frac{(b-1)A(x)}{\mu},\frac{(b-1)A(y)}{\mu},
        \frac{x-1}{x}, 1-\frac{\mu x}{A(y)B(x)}, l(y)\Bigr]$}}} $$
It's not hard to see that
\begin{align*}
 \pa_1^\infty(Z_2)& \subset\{t_3=1\}, \quad
 \pa_2^0(Z_2)\subset \{t_5=1\}, \qquad \qquad
 \pa_2^\infty(Z_2) \subset \{t_4=1\},\\
 \pa_3^\infty(Z_2)& \subset \{t_4=1\},\quad
 \pa_4^\infty(Z_2)\subset\{t_1=1\}\cup \{t_5=1\},\quad
 \pa_5^\infty(Z_2) \subset\{t_2=1\},
\end{align*}
and
\begin{align*}
\pa_1^0(Z_2)&=\Bigl[\frac{(b-1)A(y)}{\mu},
        \frac{1}{1-a}, \frac{y}{A(y)}, l(y)\Bigr]=:U''',\\
\pa_3^0(Z_2)&=\Bigl[\frac{b-1}{a\mu},\frac{(b-1)A(y)}{\mu},
       \frac{aby+1}{abA(y)}, l(y)\Bigr],\\
\pa_4^0(Z_2)&=\Bigl[\frac{(b-1)y}{B(y)},\frac{(b-1)A(y)}{\mu},
        \frac{aby+1}{abA(y)}, l(y)\Bigr],\\
\pa_5^0(Z_2)&= \Bigl[\frac{(b-1)A(x)}{\mu},\frac{(b-1)A(c)}{\mu},
        \frac{x-1}{x},\frac{y_2-x}{y_2B(x)}\Bigr]\\
\ & +\Bigl[\frac{(b-1)A(x)}{\mu},\frac{(b-1)A(y_2)}{\mu},
        \frac{x-1}{x}, \frac{c-x}{cB(x)} \Bigr].
\end{align*}
Then $U'''$ is admissible similar to $U$. By non-degeneracy
assumption it suffices to show the following cycles are
admissible:
\begin{align*}
P&:=\Bigl[\frac{(b-1)A(y)}{\mu},
       \frac{aby+1}{abA(y)}, l(y)\Bigr],\quad
Q:=\Bigl[\frac{(b-1)y}{B(y)},\frac{(b-1)A(y)}{\mu},
        \frac{aby+1}{abA(y)}, l(y)\Bigr],\\
R&:= \Bigl[\frac{(b-1)A(x)}{\mu},
        \frac{x-1}{x},\frac{y_2-x}{y_2B(x)}\Bigr],\quad
S:=\Bigl[\frac{(b-1)A(x)}{\mu},
        \frac{x-1}{x}, \frac{c-x}{cB(x)} \Bigr]. \\
\end{align*}

\begin{itemize}
\item $P$ is admissible. We have
\begin{align*}
\pa_1^0(P) &\subset \{t_3=1\},\quad
\pa_1^\infty(P)\subset\{t_2=1\},\quad
\pa_2^\infty(P)\subset\{t_3=1\},\quad
\pa_3^\infty(P)\subset\{t_1=1\},\\
\pa_2^0(P)&=\Bigl[\frac{b-1}{b},l\Big(\frac{-1}{ab}\Big)\Bigr],\quad
\pa_3^0(P)=\Bigl[\frac{(b-1)A(c)}{\mu},\frac{abc+1}{abA(c)}\Bigr]
 +\Bigl[\frac{(b-1)A(y_2)}{\mu},\frac{aby_2+1}{abA(y_2)}\Bigr].
\end{align*}
All the cycles above are admissible by \eqref{lab},
\eqref{y2nondeg}, \eqref{Ay2} and \eqref{aby2}.

\item $Q$ is admissible.  First it's easy to see that
\begin{align*}
 \pa_1^0(Q) &\subset\{t_4=1\}, \qquad
 \pa_1^\infty(Q)\subset\{t_2=1\},\quad
 \pa_2^0(Q)  \subset \{t_4=1\}, \qquad
 \pa_2^\infty(Q) \subset  \{t_3=1\},\\
 \pa_3^0(Q)&=[(1-b)/(ab-b+1),(b-1)/b,l(-1/ab)],\quad \ \,
 \pa_3^\infty(Q)\subset\{t_4=1\},\\
 \pa_4^0(Q)&=\Bigl[\frac{(b-1)c}{B(c)},\frac{(b-1)A(c)}{\mu},\frac{abc+1}{abA(c)}\Bigr]
 +\Bigl[\frac{(b-1)y_2}{B(y_2)},\frac{(b-1)A(y_2)}{\mu},\frac{aby_2+1}{abA(y_2)}\Bigr],\\
\pa_4^\infty(Q)&=[\mu/(b-1),\mu/(b-1),b].
\end{align*}
All the cycles in the above are admissible by \eqref{Ay2},
\eqref{y2nondeg}, \eqref{lab} and \eqref{aby2} and \eqref{By2}.

\item $R$ is admissible. This is because that the zeros and poles
of the three coordinate functions are all distinct.

\item $S$ is admissible. Same as $R$.
\end{itemize}

This concludes the proof that $Z_2$ is an admissible cycle.

$$\text{\boxed{\text{$Z_2(A,A)=\Bigl[\frac{(b-1)A(x)}{\mu},\frac{(b-1)A(y)}{\mu},
        \frac{x-1}{x}, \Bigl(\frac{A(y)}{y}\Bigr)
    \Bigl(1-\frac{\mu x}{A(y)B(x)}\Bigr),l(y)\Bigr]$}}} $$
We can use exactly the same proof for $Z_2$ except the following
modifications.

First, we need to look for places where we used $\{t_4=1\}$. Then
only the following needs to be re-considered:
$$\pa_3^\infty(Z_2(A,A))=\Bigl[\frac{(1-b)(1-a)}{ab-b+1},
 \frac{(b-1)A(y)}{\mu},\frac{A(y)}{y},l(y)\Bigr]$$
which can be checked to be admissible as follows: Let
\begin{equation}\label{N}
N=\Bigl[\frac{(b-1)A(y)}{\mu},\frac{A(y)}{y},l(y)\Bigr].
\end{equation}
Then
\begin{align*}
\pa_1^0(N)&\subset \{t_3=1\},\quad
 \pa_1^\infty(N)\subset\{t_2=1\},\quad
 \pa_2^0(N)\subset \{t_3=1\}, \quad
 \pa_2^\infty(N)\subset\{t_3=1\},\\
 \pa_3^0(N)&=\Bigl[\frac{(b-1)A(c)}{\mu},\frac{A(c)}{c}\Bigr]+
\Bigl[\frac{(b-1)A(y_2)}{\mu},\frac{A(y_2)}{y_2} \Bigr],\quad
 \pa_3^\infty(N)\subset\{t_1=1\}.
\end{align*} This shows that $N$, hence $\pa_3^\infty(Z_2(A,A))$, is admissible.

Second, $R$ (resp. $S$) is still admissible if we multiply the
third coordinate by $A(c)/c$ (resp. $A(y_2)/y_2$) because the
three coordinate functions of the modified cycle still have
distinct zeros and poles.

$$\text{\boxed{\text{$Z_{21}=\Bigl[\frac{(b-1)A(x)}{\mu},\frac{(b-1)A(y)}{\mu},
        \frac{x-1}{x}, \Bigl(\frac{A(y)}{y}\Bigr),l(y)\Bigr]\in C^1(F,2)\wg C^2(F,3)$}}}$$
This cycle is product of two admissible cycles
$[(b-1)A(x)/\mu,(x-1)/x]$ and $N$ given by \eqref{N}.

All the above justifies the use of Lemma 3.1(i) to get
$$Z_2(A,A)=Z_2+Z_{21}.$$

\noindent{\bf Claim 3}. Modulo admissible and negligible cycles
the two cycles $Z_3(A,A)$ and $Z_3$ are the same and both
admissible.

$$\text{\boxed{\text{$Z_3=\Bigl[A(x),A(y),\frac{abx+1}{abA(x)},
        \frac{y-x}{A(y)}, l(y)\Bigr]$}}}$$
It's not hard to see that
\begin{align*}
 \pa_1^0(Z_3)\subset &\{t_4=1\},\quad
 \pa_1^\infty(Z_3)\subset\{t_3=1\},\quad
 \pa_2^0(Z_3)\subset \{t_5=1\}, \\
 \pa_2^\infty(Z_3) \subset & \{t_4=1\},\quad
 \pa_3^\infty(Z_3) \subset \{t_4=1\}, \quad
 \pa_4^\infty(Z_3)\subset\{t_5=1\}.
\end{align*}
So we still need to show the following cycles are admissible:
\begin{align*}
\pa_3^0(Z_1)=&\Bigl[A(-1/ab),A(y),\frac{aby+1}{abA(y)},l(y)\Bigr],\quad
\pa_4^0(Z_1)=\Big[A(y),A(y),\frac{aby+1}{abA(y)}, l(y)\Big],\\
\pa_5^\infty(Z_1)= &\Bigl[A(x),\frac{\mu}{b-1},
 \frac{abx+1}{abA(x)},\frac{B(x)}{-\mu} \Bigr],\\
\pa_5^0(Z_1)=& \Bigl[A(x),A(c), \frac{abx+1}{abA(x)},
\frac{a(c-x)}{ac-a+1} \Bigr]
+\Bigl[A(x),A(y_2),\frac{abx+1}{abA(x)}, \frac{y_2-x}{A(y_2)}
\Bigr]
\end{align*}
By non-degeneracy assumption it suffices to show the following
cycles are admissible:
\begin{align*}
C'&:=\Bigl[A(y),\frac{aby+1}{abA(y)},l(y)\Bigr],\quad
C:=\Big[A(y),A(y),\frac{aby+1}{abA(y)},l(y)\Big],\quad
D:=\Bigl[A(x),\frac{abx+1}{abA(x)},\frac{B(x)}{-\mu}\Bigr],\\
E&:=\Bigl[A(x),\frac{abx+1}{abA(x)},\frac{c-x}{A(c)} \Bigr],\quad
F:=\Bigl[A(x),\frac{abx+1}{abA(x)},\frac{y_2-x}{A(y_2)}\Bigr].
\end{align*}
These are all admissible by easy computations. The only
non-obvious one identity is $y_2-(a-1)/a=A(y_2)$ which is used to
show that $F$ is admissible.

$$\text{\boxed{\text{$ Z_3(A,A)=\Bigl[\frac{(b-1)A(x)}{\mu},
 \frac{(b-1)A(y)}{\mu},\frac{abx+1}{abA(x)}, \frac{y-x}{A(y)}, l(y)\Bigr]$}}} $$
Throughout the proof that $Z_3$ is an admissible cycle we never
used the hyperplane $\{t_1=1\}$ or $\{t_2=1\}$ and moreover
$(b-1)/\mu\ne 0,\infty$ by non-degeneracy assumption. Therefore we
can use exactly the same proof to show the admissibility of
$Z_3(A,A)$.

$$\text{\boxed{\text{$Z_{31}=\Bigl[\frac{b-1}{\mu},A(y),\frac{abx+1}{abA(x)},
 \frac{y-x}{A(y)}, l(y)\Bigr]\in C^1(F,1)\wg C^2(F,4)$}}} $$
Let
$$Z'_{31}=\Bigl[A(y),\frac{abx+1}{abA(x)}, \frac{y-x}{A(y)}, l(y)\Bigr]$$
Then
\begin{align*}
 \pa_1^0(Z'_{31})&\subset \{t_4=1\},\quad
 \pa_1^\infty(Z'_{31})\subset\{t_3=1\},\quad
 \pa_2^\infty(Z'_{31})\subset \{t_3=1\},\quad
\pa_3^\infty(Z'_{31})\subset \{t_4=1\},
\end{align*}
and
\begin{align*}
 \pa_2^0(Z'_{31})&=\Bigl[A(y),\frac{aby+1}{abA(y)},l(y)\Bigr]=C',\quad
 \pa_3^0(Z'_{31})=\Bigl[A(y),\frac{aby+1}{abA(y)},l(y)\Bigr]=C', \\
\pa_4^0(Z'_{31})&=\Bigl[A(c),\frac{abx+1}{abA(x)},\frac{c-x}{A(c)}\Bigr]+
\Bigl[A(y_2),\frac{abx+1}{abA(x)}, \frac{y_2-x}{A(y_2)}\Bigr]=:E'+F',\\
\pa_4^\infty(Z'_{31})&=\Bigl[\frac{\mu}{b-1},\frac{abx+1}{abA(x)},
\frac{B(x)}{-\mu}\Bigr]=:D'.
\end{align*}
The cycles $E',F'$ and $D'$ are admissible because the coordinate
functions have different zeros and poles.

$$\text{\boxed{\text{$Z_{32}=\Bigl[A(x),\frac{b-1}{\mu},\frac{abx+1}{abA(x)},
 \frac{y-x}{A(y)}, l(y)\Bigr]\in C^1(F,1)\wg C^2(F,4)$}}} $$
Let
$$Z'_{32}=\Bigl[A(x),\frac{abx+1}{abA(x)}, \frac{y-x}{A(y)}, l(y)\Bigr] $$
Then
$$
 \pa_1^0(Z'_{32})\subset \{t_3=1\},\quad \pa_1^\infty(Z'_{32})\subset\{t_2=1\},\quad
 \pa_2^\infty(Z'_{32})\subset \{t_3=1\}, \quad
\pa_3^\infty(Z'_{32})\subset \{t_4=1\},
$$
and
\begin{align*}
 \pa_2^0(Z'_{32})&=\Bigl[A\Big(\frac{-1}{ab}\Big), \frac{aby+1}{abA(y)}, l(y)\Bigr]=:C'',\quad
 \pa_3^0(Z'_{32})=\Bigl[A(y),\frac{aby+1}{abA(y)}, l(y)\Bigr]=C',\\
\pa_4^0(Z'_{32})&=\Bigl[A(x),\frac{abx+1}{abA(x)},\frac{c-x}{A(c)}\Bigr]+
\Bigl[A(x),\frac{abx+1}{abA(x)}, \frac{y_2-x}{A(y_2)}\Bigr]=E+F,\\
\pa_4^\infty(Z'_{32})&=\Bigl[A(x),\frac{abx+1}{abA(x)},
\frac{B(x)}{-\mu}\Bigr]=D.
\end{align*}
The coordinate functions of $C''$ are all distinct because of
\eqref{y2nondeg}, \eqref{Ay2} and \eqref{aby2} so that all the
above cycles are admissible.

$$\text{\boxed{\text{$Z_{33}=\Bigl[\frac{b-1}{\mu},\frac{b-1}{\mu},
\frac{abx+1}{abA(x)}, \frac{y-x}{A(y)}, l(y)\Bigr]\in C^1(F,1)\wg
 C^1(F,1)\wg C^1(F,3)$}}} $$
Let
$$Z'_{33}=\Bigl[\frac{abx+1}{abA(x)}, \frac{y-x}{A(y)}, l(y)\Bigr]$$
Then
\begin{align*}
 \pa_1^0(Z'_{33})&=\Bigl[\frac{aby+1}{abA(y)}, l(y)\Bigr]=:C''',\quad \
 \pa_1^\infty(Z'_{33})\subset \{t_2=1\},\\
 \pa_2^0(Z'_{33})&=\Bigl[\frac{aby+1}{abA(y)}, l(y)\Bigr]=C''', \quad
\pa_2^\infty(Z'_{33})\subset \{t_3=1\}, \\
\pa_3^0(Z'_{33})&=\Bigl[\frac{abx+1}{abA(x)},\frac{c-x}{A(c)}\Bigr]+
\Bigl[\frac{abx+1}{abA(x)},\frac{y_2-x}{A(y_2)}\Bigr]=:E''+F'',\\
\pa_3^\infty(Z'_{33})&=\Bigl[\frac{abx+1}{abA(x)},
\frac{B(x)}{-\mu}\Bigr]=:D''.
\end{align*}
The admissibility of $C''',E'',F''$ and $D''$ is easy to check.

All the above justifies the use of Lemma 3.1(ii)(c)(2) to get:
$$Z_3(A,A)=Z_3+Z_{31}+Z_{32}+Z_{33}.$$

\noindent{\bf Claim 4}. Modulo admissible and negligible cycles
the two cycles $Z_4(A,A)$ and $Z_4$ are the same and both
admissible.

$$\text{\boxed{\text{$Z_4=\Bigl[\frac{(b-1)A(x)}{\mu},\frac{(b-1)A(y)}{\mu},
        \frac{abx+1}{aA(x)}, \Bigl(\frac{A(y)}{y}\Bigr)
        \Bigl(1-\frac{\mu x}{A(y)B(x)}\Bigr),l(y)\Bigr]$}}} $$
Note that
\begin{equation}\label{2ways}
\Bigl(\frac{A(y)}{y}\Bigr)\Bigl(1-\frac{\mu x}{A(y)B(x)}\Bigr)
=\frac{xB(y)+A(y)}{yB(x)}=\frac{yB(x)+A(x)}{yB(x)}.
\end{equation}
It's not hard to see that
\begin{align*}
 \pa_1^0(Z_4)&\subset\{t_4=1\},\quad \pa_2^0(Z_4)\subset \{t_5=1\}, \qquad\qquad
 \pa_2^\infty(Z_4) \subset \{t_4=1\},\\
 \pa_3^\infty(Z_4)& \subset \{t_4=1\},\quad
  \pa_4^\infty(Z_4)\subset\{t_1=1\}\cup \{t_5=1\},\quad
  \pa_5^\infty(Z_4)  \subset\{t_1=1\},
\end{align*}
and
\begin{align*}
 \pa_1^\infty(Z_4)&=\Bigl[\frac{(b-1)A(y)}{\mu},b,\frac{B(y)}{(b-1)y},l(y)\Bigr], \quad
 \pa_3^0(Z_4)=\Bigl[\frac{b-1}{b},\frac{(b-1)A(y)}{\mu},
       \frac{y-1}{y}, l(y)\Bigr],\\
 \pa_4^0(Z_4)&=\Bigl[\frac{(b-1)y}{B(y)},\frac{(b-1)A(y)}{\mu},
      \frac{y-1}{y}, l(y)\Bigr],\\
\pa_5^0(Z_4)&= \Bigl[\frac{(b-1)A(x)}{\mu},\frac{(b-1)A(c)}{\mu},
       \frac{abx+1}{aA(x)},\frac{xB(c)+A(c)}{cB(x)}\Bigr]\\
\ & +\Bigl[\frac{(b-1)A(x)}{\mu},\frac{(b-1)A(y_2)}{\mu},
       \frac{abx+1}{aA(x)}, \frac{xB(y_2)+A(y_2)}{y_2B(x)}\Bigr].
\end{align*}
By non-degeneracy assumption it suffices to show the following
cycles are admissible:
\begin{align*}
 G&:=\Bigl[\frac{(b-1)A(y)}{\mu},\frac{B(y)}{(b-1)y},l(y)\Bigr],\quad
 H:=\Bigl[\frac{(b-1)A(y)}{\mu},\frac{y-1}{y}, l(y)\Bigr],\\
 I&:=\Bigl[\frac{(b-1)y}{B(y)},\frac{(b-1)A(y)}{\mu},
      \frac{y-1}{y}, l(y)\Bigr],\quad
 J:= \Bigl[\frac{(b-1)A(x)}{\mu},\frac{abx+1}{aA(x)},\frac{xB(c)+A(c)}{cB(x)}\Bigr]\\
 K&:=\Bigl[\frac{(b-1)A(x)}{\mu},\frac{abx+1}{aA(x)},\frac{xB(y_2)+A(y_2)}{y_2B(x)}\Bigr].
\end{align*}

\begin{itemize}
\item $G$ is admissible. We have
\begin{align*}
\pa_1^0(G) &\subset \{t_3=1\},\quad
\pa_1^\infty(G)\subset\{t_2=1\},\quad
\pa_2^0(G)\subset\{t_1=1\},\quad
\pa_2^\infty(G)\subset\{t_3=1\},\\
\pa_3^0(G)&=\Bigl[\frac{(b-1)A(c)}{\mu},\frac{B(c)}{(b-1)c}\Bigr]
 +\Bigl[\frac{(b-1)A(y_2)}{\mu},\frac{B(y_2)}{(b-1)y_2}\Bigr],\quad
\pa_3^\infty(G)\subset\{t_1=1\}.
\end{align*}
It follows from \eqref{y2nondeg}, \eqref{Ay2} and \eqref{aby2}
that $\pa_3^0(G)$ is admissible.

\item $H$ is admissible. The three coordinate functions have
distinct zeros and poles because of \eqref{y2nondeg} and
\eqref{Ay2}.

\item $I$ is admissible. We have
\begin{align*}
\pa_1^0(I) &\subset \{t_4=1\},\quad
\pa_1^\infty(I)\subset\{t_2=1\},\quad
\pa_2^0(I)\subset\{t_4=1\},\quad
\pa_2^\infty(I)\subset\{t_1=1\},\\
\pa_3^0(I)&=\Bigl[\frac{b-1}{b},\frac{b-1}{a\mu},l(1)\Bigr],\quad
\pa_3^\infty(I)\subset\{t_4=1\},\quad
\pa_4^\infty(I)\subset\{t_2=1\}\\
\pa_4^0(I)&=\Bigl[\frac{(b-1)c}{B(c)},\frac{(b-1)A(c)}{\mu},\frac{c-1}{c}\Bigr]
 +\Bigl[\frac{(b-1)y_2}{B(y_2)},\frac{(b-1)A(y_2)}{\mu},\frac{y_2-1}{y_2}\Bigr].
\end{align*}
All the cycles are admissible by \eqref{Ay2}, \eqref{By2},
\eqref{l1} and \eqref{y2nondeg}.

\item $J$ and $K$ are admissible. By considering the zeros and
poles of the coordinate functions it is easy to see that the only
nontrivial thing is to check that $\pa_1^0(J)\subset \{t_3=1\}$
and $\pa_2^\infty(K)\subset\{t_3=1\}$ which follows from equation
\eqref{2ways}. For example, the zero of the 3rd coordinate of $J$
(resp. $K$) is $-A(c)/B(c)=y_2$ (resp. $-A(y_2)/B(y_2)=c$).

\end{itemize}

This concludes the proof that $Z_4$ is an admissible cycle.

$$\text{\boxed{\text{$Z_4(A,A)=\Bigl[\frac{(b-1)A(x)}{\mu},\frac{(b-1)A(y)}{\mu},
       \frac{abx+1}{abA(x)}, \Bigl(\frac{A(y)}{y}\Bigr)
    \Bigl(1-\frac{\mu x}{A(y)B(x)}\Bigr),l(y)\Bigr]$}}} $$
We can use exactly the same proof for $Z_4$ except that instead of
$I$, $J$ and $K$ we need to show $I'$, $J'$ and $K'$ are
admissible where
\begin{align*}
 I'&:=\Bigl[\frac{(b-1)y}{B(y)},\frac{(b-1)A(y)}{\mu},
      \frac{y-1}{by}, l(y)\Bigr],\\
 J'&:= \Bigl[\frac{(b-1)A(x)}{\mu},\frac{abx+1}{abA(x)},\frac{xB(c)+A(c)}{cB(x)}\Bigr]\\
 K'&:=\Bigl[\frac{(b-1)A(x)}{\mu},\frac{abx+1}{abA(x)},\frac{xB(y_2)+A(y_2)}{y_2B(x)}\Bigr].
\end{align*}
Exactly the same proofs are valid because in the proof of $I$ we
didn't use the hyperplane $\{t_3=1\}$ while for $J$ and $K$ we
didn't use $\{t_2=1\}$.

$$\text{\boxed{\text{$Z_{41}=\Bigl[\frac{(b-1)A(x)}{\mu},\frac{(b-1)A(y)}{\mu},
        \frac{1}{b}, \Bigl(\frac{A(y)}{y}\Bigr)
        \Bigl(1-\frac{\mu x}{A(y)B(x)}\Bigr),l(y)\Bigr]\in C^1(F,1)\wg C^2(F,4)$}}}$$
The same argument for $Z_4(A,A)$ goes through without any problem.

All the above justifies the use of Lemma 3.1(ii)(a) to get
$$Z_4(A,A)=Z_4+Z_{41}.$$

\medskip
\noindent {\bf Step} (4). Decomposition of $\rho_x Z_2(A,A)+\rho_y
Z_4(A,A)=X_1-X_2$.

We need to show
\begin{equation*}
    X_1=X_{11}+X_{12}, \quad X_2=X_{21}+X_{22},
\end{equation*}
where
\begin{align*}
X_{11}=&\Bigl[\frac{(b-1)x}{B(x)},\frac{A(y)}{-\mu y},
    \frac{abx+1}{aA(x)},\frac{y-x}{A(y)},l(y)\Bigr]\\
X_{12}=&\Bigl[\frac{A(x)}{-\mu x},\frac{(b-1)y}{B(y)},
\frac{abx+1}{aA(x)},
    \frac{\mu(x-y)}{A(y)B(x)},l(y)\Bigr],\\
X_{21}=&\Bigl[\frac{B(x)}{(b-1)x},(1-b)y,\frac{abx+1}{aA(x)},
    \frac{y-x}{A(y)},l(y)\Bigr]\\
X_{22}=&\Bigl[(1-b)x,\frac{B(y)}{(b-1)y}, \frac{abx+1}{aA(x)},
    \frac{\mu(x-y)}{A(y)B(x)},l(y)\Bigr],
\end{align*}
are admissible.

$$\text{\boxed{\text{$X_{11}=\Bigl[\frac{(b-1)x}{B(x)},\frac{A(y)}{-\mu y},
    \frac{abx+1}{aA(x)},\frac{y-x}{A(y)},l(y)\Bigr]$}}} $$
We have
\begin{alignat*}{5}
 \pa_1^\infty(X_{11})& \subset\{t_3=1\},\quad&
 \pa_2^0(X_{11})&\subset \{t_5=1\}, \quad&
 \pa_2^\infty(X_{11})& \subset \{t_5=1\},\\
 \pa_3^\infty(X_{11})& \subset \{t_4=1\},\quad&
 \pa_4^\infty(X_{11})& \subset\{t_5=1\},\quad&
 \pa_5^\infty(X_{11})& \subset\{t_2=1\},
\end{alignat*}
and
\begin{align*}
 \pa_1^0(X_{11})&=\Bigl[\frac{A(y)}{-\mu y},\frac{1}{1-a},\frac{y}{A(y)},l(y)\Bigr]=:U^{(4)},\\
 \pa_3^0(X_{11})&=\Bigl[\frac{b-1}{a\mu},\frac{A(y)}{-\mu y},
       \frac{aby+1}{abA(y)}, l(y)\Bigr]=:P',\\
 \pa_4^0(X_{11})&=\Bigl[\frac{(b-1)y}{B(y)},\frac{A(y)}{-\mu y},
    \frac{aby+1}{aA(y)},l(y)\Bigr]=:Q',\\
 \pa_5^0(X_{11})&=\Bigl[\frac{(b-1)x}{B(x)},\frac{A(c)}{-c\mu},
    \frac{abx+1}{aA(x)},\frac{c-x}{A(c)}\Bigr]\\
\ & +\Bigl[\frac{(b-1)x}{B(x)},\frac{A(y_2)}{-y_2\mu},
    \frac{abx+1}{aA(x)},\frac{y_2-x}{A(y_2)}\Bigr]=:M_1+M_2.
\end{align*}
All these cycles are admissible by arguments similarly to those
for $U$, $P$, and $Q$. For $M_i$ ($i=1,2$) we can see that the
coordinate functions have distinct zeros and poles.

$$\text{\boxed{\text{$X_{12}=\Bigl[\frac{A(x)}{-\mu x},\frac{(b-1)y}{B(y)},
 \frac{abx+1}{aA(x)},\frac{\mu(x-y)}{A(y)B(x)},l(y)\Bigr]$}}}$$
We have
\begin{alignat*}{5}
 \pa_1^0(X_{12})& \subset\{t_4=1\},\quad&
 \pa_2^0(X_{12})&\subset \{t_5=1\}, \quad&
 \pa_2^\infty(X_{12})& \subset \{t_4=1\},\\
 \pa_3^\infty(X_{12})& \subset \{t_4=1\},\quad&
 \pa_4^\infty(X_{12})& \subset\{t_5=1\}\cup\{t_3=1\},\quad&
 \pa_5^\infty(X_{12})& \subset\{t_4=1\},
\end{alignat*}
and
\begin{align*}
 \pa_1^\infty(X_{12})&=\Bigl[\frac{(b-1)y}{B(y)},\frac{1}{1-a},\frac{-\mu y}{A(y)},l(y)\Bigr]=:Q''
,\\
 \pa_3^0(X_{12})&=\Bigl[a,\frac{(b-1)y}{B(y)},
       \frac{aby+1}{aA(y)}, l(y)\Bigr]=:Q''',\\
 \pa_4^0(X_{12})&=\Bigl[\frac{A(y)}{-\mu y},\frac{(b-1)y}{B(y)},
    \frac{aby+1}{aA(y)},l(y)\Bigr]=-Q',\\
 \pa_5^0(X_{12})&=\Bigl[\frac{A(x)}{-\mu x},\frac{(b-1)c}{B(c)},
    \frac{abx+1}{aA(x)},\frac{\mu(x-c)}{A(c)B(x)}\Bigr]\\
\ & +\Bigl[\frac{A(x)}{-\mu x},\frac{(b-1)y_2}{B(y_2)},
    \frac{abx+1}{aA(x)},\frac{\mu(x-y_2)}{A(y_2)B(x)}\Bigr]=:N_1+N_2.
\end{align*}
Both $Q''$ and $Q'''$ are admissible by argument similarly to that
for $Q$. For $N_i$ ($i=1,2$) we can see that the coordinate
functions have distinct zeros and poles except when $A(x)=0$ which
implies that $t_4=1.$

$$\text{\boxed{\text{$X_{21}=\Bigl[\frac{B(x)}{(b-1)x},(1-b)y,\frac{abx+1}{aA(x)},
    \frac{y-x}{A(y)},l(y)\Bigr]$}}} $$
We have
\begin{alignat*}{5}
 \pa_1^0(X_{21})& \subset\{t_3=1\},\quad&
 \pa_2^0(X_{21})&\subset \{t_5=1\}, \quad&
 \pa_2^\infty(X_{21})& \subset \{t_4=1\},\\
 \pa_3^\infty(X_{21})& \subset \{t_4=1\},\quad&
 \pa_4^\infty(X_{21})& \subset\{t_5=1\},\quad&
 \pa_5^\infty(X_{21})& \subset\{t_2=1\},
\end{alignat*}
and
\begin{align*}
 \pa_1^\infty(X_{21})&=\Bigl[(1-b)y,\frac{1}{1-a},\frac{y}{A(y)},l(y)\Bigr]=:U^{(5)}
,\\
 \pa_3^0(X_{21})&=\Bigl[\frac{a\mu}{b-1},(1-b)y,
       \frac{aby+1}{abA(y)}, l(y)\Bigr]=:P'',\\
 \pa_4^0(X_{21})&=\Bigl[\frac{B(y)}{(b-1)y},(1-b)y,
    \frac{aby+1}{aA(y)},l(y)\Bigr]=:Q^{(4)},\\
 \pa_5^0(X_{21})&=\Bigl[\frac{B(x)}{(b-1)x},(1-b)c,
    \frac{abx+1}{aA(x)},\frac{c-x}{A(c)}\Bigr]\\
\ & +\Bigl[\frac{B(x)}{(b-1)x},(1-b)y_2,
    \frac{abx+1}{aA(x)},\frac{y_2-x}{A(y_2)}\Bigr]=:O_1+O_2.
\end{align*}
All these cycles are admissible by arguments similarly to those
for $U$, $P$, and $Q$. For $O_i$ ($i=1,2$) we can see that the
coordinate functions have distinct zeros and poles.

$$\text{\boxed{\text{$X_{22}=\Bigl[(1-b)x,\frac{B(y)}{(b-1)y}, \frac{abx+1}{aA(x)},
    \frac{\mu(x-y)}{A(y)B(x)},l(y)\Bigr]$}}} $$
We have
\begin{alignat*}{3}
 \pa_2^0(X_{22})&\subset \{t_4=1\}, \quad
 \pa_2^\infty(X_{22}) \subset \{t_5=1\},\quad&
 \pa_3^\infty(X_{22})\subset \{t_4=1\},\\
 \pa_4^\infty(X_{22})& \subset\{t_5=1\}\cup   \{t_3=1\},\quad&
 \pa_5^\infty(X_{22})\subset\{t_4=1\},
\end{alignat*}
and
\begin{align*}
 \pa_1^0(X_{22})&=\Bigl[\frac{B(y)}{(b-1)y},\frac{1}{1-a},\frac{-\mu y}{A(y)},l(y)\Bigr]=-Q'',\\
 \pa_1^\infty(X_{22})&=\Bigl[\frac{B(y)}{(b-1)y},b,\frac{\mu}{(b-1)A(y)},l(y)\Bigr]=:Q^{(4)},\\
 \pa_3^0(X_{22})&=\Bigl[\frac{b-1}{ab},\frac{B(y)}{(b-1)y},
       \frac{aby+1}{aA(y)}, l(y)\Bigr]=-Q''',\\
 \pa_4^0(X_{22})&=\Bigl[(1-b)y,\frac{B(y)}{(b-1)y},
    \frac{aby+1}{aA(y)},l(y)\Bigr]=:Q^{(5)},\\
 \pa_5^0(X_{22})&=\Bigl[(1-b)x,\frac{B(c)}{(b-1)c},
    \frac{abx+1}{aA(x)},\frac{\mu(x-c)}{A(c)B(x)}\Bigr]\\
\ & +\Bigl[(1-b)x,\frac{B(y_2)}{(b-1)y_2},
    \frac{abx+1}{aA(x)},\frac{\mu(x-y_2)}{A(y_2)B(x)}\Bigr]=:P_1+P_2.
\end{align*}
All these cycles are admissible by argument similarly to that for
$Q$. For $P_i$ ($i=1,2$) we can consider the coordinate functions
and see that they have all distinct zeros and poles.

\medskip
\noindent {\bf Step} (5). Computation of  $X_1$.

\medskip

Set
$$\tZ(f_1,f_2)=\Bigl[f_1,f_2,\frac{abx+1}{aA(x)},
    \frac{\mu(x-y)}{A(y)B(x)},l(y)\Bigr].$$
We want to show that by throwing away the appropriate admissible
and negligible cycle we have
$$X_1=\tZ\Bigl(\frac{(b-1)f}{B},\frac{A}{-\mu f}\Bigr)
    +\tZ\Bigl(\frac{A}{-\mu f},\frac{(b-1)f}{B}\Bigr).$$
For this step we need to use Lemma 3.1(i) to get
\begin{align*}
X_{13}:=&\Bigl[\frac{(b-1)x}{B(x)},\frac{A(y)}{-\mu y},
\frac{abx+1}{aA(x)},
    \frac{\mu(x-y)}{A(y)B(x)},l(y)\Bigr]\\
=&\Bigl[\frac{(b-1)x}{B(x)},\frac{A(y)}{-\mu y},
\frac{abx+1}{aA(x)},
    \frac{y-x}{A(y)},l(y)\Bigr]\\
+&\Bigl[\frac{(b-1)x}{B(x)},\frac{A(y)}{-\mu y},
 \frac{abx+1}{aA(x)},\frac{-\mu}{B(x)},l(y)\Bigr]=:X_{11}+X_{14}.
\end{align*}
So it suffices to show that both $X_{13}$ and $X_{14}$ are
admissible. It's obvious that $X_{14}$ is the product of two
admissible cycles
$$\Bigl[\frac{(b-1)x}{B(x)},\frac{abx+1}{aA(x)},\frac{-\mu}{B(x)}\Bigr],\quad
 \Bigl[\frac{A(y)}{-\mu y},l(y)\Bigr]$$
while the admissibility of $X_{13}$ can be shown by the same
argument as that for $X_{11}$ except for the last step
\begin{align*}
 \pa_5^0(X_{13})&=\Bigl[\frac{(b-1)x}{B(x)},\frac{aA(c)}{(ab-b+1)c},
    \frac{abx+1}{aA(x)},\frac{\mu(x-c)}{A(c)B(x)}\Bigr]\\
\ & +\Bigl[\frac{(b-1)x}{B(x)},\frac{aA(y_2)}{(ab-b+1)y_2},
    \frac{abx+1}{aA(x)},\frac{\mu(x-y_2)}{A(y_2)B(x)}\Bigr]=:R_1+R_2.
\end{align*}
By consideration of the zeros and poles of the coordinate
functions we can show that $R_i$ ($i=1,2$) is admissible because
$B(x)=0$ implies $t_3=1$.

Next we want to show
$$Z_3(F,F)=\tZ(F,F) \text{ for }F=\frac{A}{f},\
    \frac{f}{B},\ \frac{A}{B}$$
where
\begin{align}\label{Z3fB}
Z_3\Bigl(\frac{f}{B},\frac{f}{B}\Bigr)
    &:=\Bigl[\frac{(b-1)x}{B(x)},\frac{(1-b)y}{B(y)},
    \frac{abx+1}{aA(x)},\frac{y-x}{yB(x)},l(y)\Bigr],\\
Z_3\Bigl(\frac{A}{f},\frac{A}{f}\Bigr)
    &:=\Bigl[\frac{A(x)}{-\mu x},\frac{A(y)}{-\mu y},
    \frac{abx+1}{aA(x)}, \frac{y-x}{A(y)},l(y)\Bigr],\label{Z3Af}\\
Z_3\Bigl(\frac{A}{B},\frac{A}{B}\Bigr)
    &:=\Bigl[\frac{A(x)}{B(x)},\frac{A(y)}{B(y)},\frac{abx+1}{aA(x)},
    \frac{\mu(x-y)}{A(y)B(x)},l(y)\Bigr].\label{Z3AB}
\end{align}
To prove these it suffices to show the following: First,
\begin{align*}
 \ &\Bigl[\frac{(b-1)x}{B(x)},\frac{(1-b)y}{B(y)},
    \frac{abx+1}{aA(x)},\frac{-\mu y}{A(y)},l(y)\Bigr]\in C^1(F,2)\wg C^2(F,3),\\
 \ &\Bigl[\frac{A(x)}{-\mu x},\frac{A(y)}{-\mu y},
    \frac{abx+1}{aA(x)},\frac{-a\mu}{B(x)},l(y)\Bigr]\in C^1(F,2)\wg C^2(F,3)
\end{align*}
are both admissible and negligible which is not hard to see. Then
all the following are admissible and negligible:
\begin{align*}
\ &\Bigl[b-1,\frac{(1-b)y}{B(y)},
    \frac{abx+1}{aA(x)},\frac{y-x}{yB(x)},l(y)\Bigr]\in C^1(F,1)\wg C^2(F,4),\\
\ &\Bigl[\frac{(b-1)x}{B(x)},b-1,
    \frac{abx+1}{aA(x)},\frac{y-x}{yB(x)},l(y)\Bigr]\in C^1(F,1)\wg C^2(F,4),\\
\ &\Bigl[\frac{-1}{\mu},\frac{A(y)}{-\mu y},
    \frac{abx+1}{aA(x)}, \frac{y-x}{A(y)},l(y)\Bigr]\in C^1(F,1)\wg C^2(F,4),\\
\ &\Bigl[\frac{A(x)}{-\mu x},\frac{-1}{\mu},
    \frac{abx+1}{aA(x)}, \frac{y-x}{A(y)},l(y)\Bigr]\in C^1(F,1)\wg C^2(F,4),\\
\ &\Bigl[\mu,\mu,\frac{abx+1}{aA(x)},
\frac{y-x}{A(y)},l(y)\Bigr]\in C^1(F,1)\wg C^1(F,1)\wg C^1(F,3)\
\text{for any }\mu\ne 0.
\end{align*}
We only need to note that $B(x)=0$ implies $t_3=1$ and that
$yA(y)=0$ implies $t_5=1$ for all the above cycles, and $B(y)=0$
implies that $t_4=1$ for the first two cycles.

\medskip
\noindent {\bf Step} (6). Decomposition of $X_2=Y_1+Y_2+Y_3+Y_4$.

\medskip
Put
$$v(x)=\frac{abx+1}{aA(x)},
\quad  l_1(y)=1-\frac{y}{c}, \quad l_2(y)=\frac{y_2-y}{y_2B(y)},
$$ which satisfies
$$l_1(y)l_2(y)=l(y)=1-\frac{k(c)}{k(y)},\quad l_1(0)=l_2(0)=1.$$
Then it follows from Lemma~3.2(i) that
\begin{equation}\label{X2}
X_2=Y_1+Y_2+Y_3+Y_4
\end{equation}
where all of the cycles $$\aligned Y_1=&\Bigl[\frac{B(x)}{(b-1)x},
(b-1)y,\frac{abx+1}{aA(x)},
    \frac{y-x}{A(y)},l_1(y)\Bigr],\\
Y_2=&\Bigl[(b-1)x,\frac{B(y)}{(b-1)y},\frac{abx+1}{aA(x)},
    \frac{\mu(x-y)}{A(y)B(x)},l_1(y)\Bigr],\\
Y_3=&\Bigl[\frac{(b-1)x}{B(x)}, (b-1)y,\frac{abx+1}{aA(x)},
    \frac{y-x}{A(y)},l_2(y)\Bigr],\\
Y_4=&\Bigl[(b-1)x,\frac{(b-1)y}{B(y)},\frac{abx+1}{aA(x)},
    \frac{\mu(x-y)}{A(y)B(x)},l_2(y)\Bigr]
\endaligned$$
are admissible. This breakup is the key step in the whole paper.

$$\text{\boxed{\text{$Y_1=\Bigl[\frac{B(x)}{(b-1)x},(b-1)y,\frac{abx+1}{aA(x)},
    \frac{y-x}{A(y)},l_1(y)\Bigr],\\
$}}} $$ We have
\begin{alignat*}{2}
 \pa_1^0(Y_1)&\subset \{t_3=1\}, \quad
 \pa_2^0(Y_1)\subset \{t_5=1\}, \quad
 \pa_2^\infty(Y_1) \subset \{t_4=1\},\\
 \pa_3^\infty(Y_1)&\subset \{t_4=1\},\quad
 \pa_5^\infty(Y_1)\subset\{t_4=1\},
\end{alignat*}
and
\begin{align*}
 \pa_1^\infty(Y_1)&=\Bigl[(b-1)y,\frac{1}{1-a},\frac{y}{A(y)},1-\frac{y}{c}\Bigr],\\
 \pa_3^0(Y_1)&=\Bigl[\frac{a\mu}{b-1},(b-1)y,
       \frac{aby+1}{abA(y)},1-\frac{y}{c}\Bigr],\\
 \pa_4^0(Y_1)&=\Bigl[\frac{B(y)}{(b-1)y},(b-1)y,
    \frac{aby+1}{aA(y)},1-\frac{y}{c}\Bigr],\\
 \pa_4^\infty(Y_1)&=\Bigl[\frac{B(x)}{(b-1)x},\frac{(b-1)(a-1)}{a},
    \frac{abx+1}{aA(x)},\frac{ac-a+1}{ac}\Bigr],\\
 \pa_5^0(Y_1)&=\Bigl[\frac{B(x)}{(b-1)x},(b-1)c,
    \frac{abx+1}{aA(x)},\frac{c-x}{aA(c)}\Bigr].
\end{align*}
All these cycles are clearly admissible.

$$\text{\boxed{\text{$Y_2=\Bigl[(b-1)x,\frac{B(y)}{(b-1)y},\frac{abx+1}{aA(x)},
    \frac{\mu(x-y)}{A(y)B(x)},l_1(y)\Bigr]$}}} $$
We have
$$ \pa_2^0(Y_2)\subset \{t_4=1\}, \quad
 \pa_2^\infty(Y_2) \subset \{t_5=1\},\quad
 \pa_3^\infty(Y_2)\subset \{t_4=1\},\quad
 \pa_5^\infty(Y_2)\subset\{t_2=1\},
$$
and
\begin{align*}
 \pa_1^0(Y_2)&=\Bigl[\frac{B(y)}{(b-1)y},\frac{1}{1-a},
    \frac{-\mu y}{A(y)},1-\frac{y}{c}\Bigr],\\
 \pa_1^\infty(Y_2)&=\Bigl[\frac{B(y)}{(b-1)y},b,
    \frac{\mu}{(b-1)A(y)},1-\frac{y}{c}\Bigr],\\
 \pa_3^0(Y_2)&=\Bigl[\frac{1-b}{ab},\frac{B(y)}{(b-1)y},
       \frac{aby+1}{aA(y)},1-\frac{y}{c}\Bigr],\\
 \pa_4^0(Y_2)&=\Bigl[(b-1)y,\frac{B(y)}{(b-1)y},
    \frac{aby+1}{aA(y)},\frac{c-y}{c}\Bigr],\\
 \pa_4^\infty(Y_2)&=\Bigl[(b-1)x,\frac{ab-b+1}{(b-1)(a-1)},
    \frac{abx+1}{aA(x)},\frac{ac-a+1}{ac}\Bigr],\\
 \pa_5^0(Y_2)&=\Bigl[(b-1)x,\frac{B(c)}{(b-1)c},\frac{abx+1}{aA(x)},
    \frac{\mu(x-c)}{A(c)B(x)}\Bigr].
\end{align*}
All these cycles are clearly admissible.

$$\text{\boxed{\text{$Y_3=\Bigl[\frac{(b-1)x}{B(x)}, (b-1)y,\frac{abx+1}{aA(x)},
    \frac{y-x}{A(y)},l_2(y)\Bigr],\\
$}}} $$ We have
\begin{alignat*}{2}
 \pa_1^\infty(Y_3)&\subset \{t_3=1\}, \quad
 \pa_2^0(Y_3)\subset \{t_5=1\}, \quad
 \pa_2^\infty(Y_3) \subset \{t_4=1\},\\
 \pa_3^\infty(Y_3)&\subset \{t_4=1\},\quad
 \pa_5^\infty(Y_3)\subset\{t_4=1\},
\end{alignat*}
and
\begin{align*}
 \pa_1^0(Y_3)&=\Bigl[(b-1)y,\frac{1}{1-a},
    \frac{y}{A(y)},\frac{y_2-y}{y_2B(y)}\Bigr],\\
 \pa_3^0(Y_3)&=\Bigl[\frac{b-1}{a\mu},(b-1)y,
       \frac{aby+1}{abA(y)},\frac{y_2-y}{y_2B(y)}\Bigr],\\
 \pa_4^0(Y_3)&=\Bigl[\frac{(b-1)y}{B(y)},(b-1)y,
    \frac{aby+1}{aA(y)},\frac{y_2-y}{y_2B(y)}\Bigr],\\
 \pa_4^\infty(Y_3)&=\Bigl[\frac{(b-1)x}{B(x)},\frac{(b-1)(a-1)}{a},
    \frac{abx+1}{aA(x)},\frac{ac}{ac-a+1}\Bigr],\\
 \pa_5^0(Y_3)&=\Bigl[\frac{(b-1)x}{B(x)},(b-1)y_2,
    \frac{abx+1}{aA(x)},\frac{y_2-x}{A(y_2)}\Bigr].
\end{align*}
All these cycles are clearly admissible.

$$\text{\boxed{\text{$Y_4=\Bigl[(b-1)x,\frac{(b-1)y}{B(y)},\frac{abx+1}{aA(x)},
    \frac{\mu(x-y)}{A(y)B(x)},l_2(y)\Bigr]$}}} $$
We have
$$ \pa_2^0(Y_4)\subset \{t_5=1\}, \quad
 \pa_2^\infty(Y_4) \subset \{t_4=1\},\quad
 \pa_3^\infty(Y_4)\subset \{t_4=1\},\quad
 \pa_5^\infty(Y_4)\subset\{t_2=1\},
$$
and
\begin{align*}
 \pa_1^0(Y_4)&=\Bigl[\frac{(b-1)y}{B(y)},\frac{1}{1-a},
    \frac{-\mu y}{A(y)},\frac{y_2-y}{y_2B(y)}\Bigr],\\
 \pa_1^\infty(Y_4)&=\Bigl[\frac{(b-1)y}{B(y)},b,
    \frac{\mu}{(b-1)A(y)},\frac{y_2-y}{y_2B(y)}\Bigr],\\
 \pa_3^0(Y_4)&=\Bigl[\frac{1-b}{ab},\frac{(b-1)y}{B(y)},
       \frac{aby+1}{aA(y)},\frac{y_2-y}{y_2B(y)}\Bigr],\\
 \pa_4^0(Y_4)&=\Bigl[(b-1)y,\frac{(b-1)y}{B(y)},
    \frac{aby+1}{aA(y)},\frac{y_2-y}{y_2B(y)}\Bigr],\\
 \pa_4^\infty(Y_4)&=\Bigl[(b-1)x,\frac{(b-1)(a-1)}{ab-b+1},
    \frac{abx+1}{aA(x)},\frac{ac}{ac-a+1}\Bigr],\\
 \pa_5^0(Y_4)&=\Bigl[(b-1)x,\frac{(b-1)y_2}{B(y_2)},\frac{abx+1}{aA(x)},
    \frac{\mu(x-y_2)}{A(y_2)B(x)}\Bigr].
\end{align*}
All these cycles are clearly admissible.

\medskip
\noindent {\bf Step} (7). Computation of $Y_1+Y_2$.

\medskip

Set
$$\ga=\frac{bc-c}{bc-c+1}, \qquad \gd= \frac{1}{b},$$
and
$${\setlength\arraycolsep{1pt}
\begin{array}{rlrlrl}
v(x)&={\displaystyle \frac{abx+1}{aA(x)}}, & \quad
g(x)&={\displaystyle \frac{B(x)}{(b-1)x}}, &\quad h(x)&=(b-1)x,
\\
\quad p_4(x,y)&={\displaystyle \frac{\mu(x-y)}{A(y)B(x)}}, & \quad
q_4(x,y)&={\displaystyle \frac{y-x}{A(y)}}, & \quad \quad
s_4(x,y)&={\displaystyle \frac{(b-1)(y-x)}{B(y)}},
\\
r_4(x,y)&={\displaystyle \frac{(b-1)(y-x)}{xB(y)}}, & \quad
w_4(x,y)&={\displaystyle \frac{y-x}{B(x)(y-1)}}.
\end{array}}$$
such that $\ga l_1\big(1/(1-b)\big)=\gd v(\infty)=1$. By Lemma
3.1(ii)(1) we get
\begin{align*}
2[gh, gh,\gd v, q_4,\ga l_1]
=&[gh, gh,\gd v, q_4,\ga l_1]+[gh, gh,\gd v, s_4,\ga l_1]\\
 =&[g, gh,\gd v, q_4, \ga l_1]+[h, gh,\gd v, q_4, \ga l_1]\\
    \ +&[gh, g,\gd v, s_4, \ga l_1]+[gh, h,\gd v, s_4, \ga l_1]
\end{align*}
are all admissible. Then applying Lemma~3.1 and Lemma~3.2
repeatedly we get
\begin{align*}
 \ &[g, gh,\gd v, q_4, \ga l_1]+[gh, g,\gd v, s_4, \ga l_1]\\
 =&[g,gh, \gd v, q_4, \ga l_1]+[gh, g, \gd v, r_4, \ga l_1]\\
 =&[g,gh, v, q_4, \ga l_1]+[gh, g,  v, r_4, \ga l_1]\\
 =&[g,gh, v, q_4, \ga l_1]+[gh, g,  v, w_4, \ga l_1]\\
 =&[g, gh,v, q_4,  l_1]+[gh, g, v, w_4, l_1]\\
 =&[g, gh,v, q_4,  l_1]+[gh, g, v, p_4, l_1]\\
 =&[g,h,v, q_4, l_1]+[h, g, v, p_4, l_1]
 +[g, g, v,q_4,l_1]+[g,g, v, p_4,l_1]\\
 =&[g,h,v, q_4, l_1]+[h, g, v, p_4, l_1]
 +2[g, g, v, p_4,l_1].
\end{align*}
Again applying Lemma~3.2(i) and (ii) and Lemma~3.1(ii), we have
\begin{align*}
 \ &[h, gh,\gd v,q_4,\ga l_1]+[gh, h,\gd v, s_4,\ga l_1]\\
 =&[h, gh,\gd v,q_4, l_1]+[gh, h,\gd v, s_4, l_1]\\
 =&[h, gh,\gd v,q_4, l_1]+[gh, h,\gd v, q_4, l_1]\\
 =&[h, g, \gd v,q_4, l_1]+[g, h,\gd v, q_4, l_1]
 +2[h, h,\gd v,q_4, l_1]\\
 =&[h, g, v, q_4, l_1]+[g, h, v, q_4, l_1]
 +2[h, h,\gd v,q_4, l_1]\\
 =&[h, g, v, p_4, l_1]+[g, h, v, q_4, l_1]
 +2[h, h,\gd v,q_4, l_1].
\end{align*}

Let's prove the admissibility of all the cycles appearing in these
equations.

\begin{itemize}
\item $[gh, gh,\gd v, q_4,\ga l_1]$ is admissible.

$$\text{\boxed{\text{$[B,B]:=[gh, gh,\gd v, q_4,\ga l_1]
 =\Bigl[B(x),B(y), \frac{abx+1}{abA(x)},\frac{y-x}{A(y)},\ga l_1(y)\Bigr]$}}} $$
We have
\begin{align*}
 \pa_1^\infty[B,B]&\subset \{t_3=1\}, \quad
 \pa_2^0[B,B]\subset \{t_5=1\}, \quad
 \pa_2^\infty[B,B] \subset \{t_4=1\},\\
 \pa_3^\infty[B,B]&\subset \{t_4=1\},\quad
 \pa_5^\infty[B,B]\subset\{t_4=1\},
\end{align*}
and
\begin{align*}
 \pa_1^0[B,B]&=\Bigl[B(y),\frac{1}{b},
    \frac{B(y)}{(b-1)A(y)},\ga l_1(y) \Bigr],\\
 \pa_3^0[B,B]&=\Bigl[\frac{\mu}{b},B(y),
       \frac{aby+1}{abA(y)},\ga l_1(y) \Bigr],\\
 \pa_4^0[B,B]&=\Bigl[B(y),B(y),
    \frac{aby+1}{abA(y)},\ga l_1(y) \Bigr],\\
 \pa_4^\infty[B,B]&=\Bigl[B(x),\mu,
    \frac{abx+1}{abA(x)}, \frac{(b-1)(ac-a+1)}{a(bc-c+1)}\Bigr],\\
 \pa_5^0[B,B]&=\Bigl[B(x),B(c),\frac{abx+1}{abA(x)},
    \frac{c-x}{A(c)}\Bigr].
\end{align*}
All these cycles are clearly admissible by our choice of $\ga$ and
$\gd$.

\item $[g, gh, \cdots ]$ are admissible.

$$\text{\boxed{\text{$[B/f,B]:=[g, gh,\gd v, q_4, \ga l_1]
 =\Bigl[\frac{B(x)}{(b-1)x},B(y), \frac{abx+1}{abA(x)},\frac{y-x}{A(y)},\ga l_1(y)\Bigr]$}}} $$
We have
\begin{align*}
 \pa_2^0[B/f,B]&\subset \{t_5=1\}, \quad
 \pa_2^\infty[B/f,B] \subset \{t_4=1\},\\
 \pa_3^\infty[B/f,B]&\subset \{t_4=1\},\quad
 \pa_5^\infty[B/f,B]\subset\{t_4=1\},
\end{align*}
and
\begin{align*}
 \pa_1^0[B/f,B]&=\Bigl[B(y),\frac{1}{b},
    \frac{B(y)}{(b-1)A(y)},\ga l_1(y) \Bigr],\\
 \pa_1^\infty[B/f,B]&=\Bigl[B(y),\frac{1}{b(1-a)},
    \frac{y}{A(y)},\ga l_1(y) \Bigr],\\
 \pa_3^0[B/f,B]&=\Bigl[\frac{a\mu}{b-1},B(y),
       \frac{aby+1}{abA(y)},\ga l_1(y) \Bigr],\\
 \pa_4^0[B/f,B]&=\Bigl[\frac{B(y)}{(b-1)y},B(y),
    \frac{aby+1}{abA(y)},\ga l_1(y) \Bigr],\\
 \pa_4^\infty[B/f,B]&=\Bigl[\frac{B(x)}{(b-1)x},\mu,
    \frac{abx+1}{abA(x)}, \frac{(b-1)(ac-a+1)}{a(bc-c+1)}\Bigr],\\
 \pa_5^0[B/f,B]&=\Bigl[\frac{B(x)}{(b-1)x},B(c),\frac{abx+1}{abA(x)},
    \frac{c-x}{A(c)}\Bigr].
\end{align*}
All these cycles are clearly admissible by our choice of $\ga$.
Then by similar argument we can see that
$$[g, gh,v, q_4, \ga l_1], \quad \text{and}\quad [g, gh,\gd, q_4, \ga l_1]$$
are both admissible because we didn't use $t_3=1$ in the above.

$$\text{\boxed{\text{$[B/f,B]_1:=[g, gh, v, q_4, l_1]
 =\Bigl[\frac{B(x)}{(b-1)x},B(y), \frac{abx+1}{aA(x)},\frac{y-x}{A(y)}, l_1(y)\Bigr]$}}} $$
The same proof as above except whenever we used $\ga
l_1(1/(1-b))=1$ (namely $t_5=1$) before we have to use
$v(1/(1-b))=1$ now.

\item $[gh, g,\cdots]$ are admissible.

$$\text{\boxed{\text{$[B,B/f]:=[gh, g,\gd v, s_4, \ga l_1]
 =\Bigl[B(x),\frac{B(y)}{(b-1)y}, \frac{abx+1}{abA(x)},\frac{(b-1)(y-x)}{B(y)},\ga l_1(y)\Bigr]$}}} $$
We have
\begin{alignat*}{3}
 \pa_1^0[B,B/f]& \subset \{t_4=1\},\quad&
 \pa_1^\infty[B,B/f]&\subset \{t_3=1\},\quad
  \pa_2^0[B,B/f]&\subset \{t_5=1\}, \\
\pa_4^\infty[B,B/f]&\subset\{t_5=1\}, \quad &
 \pa_5^\infty[B,B/f]&\subset\{t_2=1\},
\end{alignat*}
and
\begin{align*}
 \pa_2^\infty[B,B/f]&=\Bigl[B(x),\frac{abx+1}{abA(x)},(1-b)x,\ga \Bigr],\\
 \pa_3^0[B,B/f]&=\Bigl[\frac{\mu}{b},\frac{B(y)}{(b-1)y},
       \frac{(b-1)(aby+1)}{abB(y)},\ga l_1(y) \Bigr],\\
 \pa_3^\infty[B,B/f]&=\Bigl[-\mu,\frac{B(y)}{(b-1)y},
       \frac{(b-1)A(y)}{B(y)},\ga l_1(y) \Bigr],\\
 \pa_4^0[B,B/f]&=\Bigl[B(y),\frac{B(y)}{(b-1)y},
   \frac{aby+1}{abA(y)},\ga l_1(y) \Bigr],\\
  \pa_5^0[B,B/f]&=\Bigl[B(x),\frac{B(c)}{(b-1)c},\frac{abx+1}{abA(x)},
   \frac{(b-1)(c-x)}{B(c)}\Bigr].
\end{align*}
All these cycles are clearly admissible by our choice of $\ga$.

$$\text{\boxed{\text{$[B,B/f]':=[gh, g,\gd v, s_4/r_4, \ga l_1]=\Bigl[B(x),\frac{B(y)}{(b-1)y},
 \frac{abx+1}{abA(x)},x,\ga l_1(y)\Bigr]$}}} $$
This a product of two admissible cycles.

$$\text{\boxed{\text{$[B,B/f]_1:=[gh, g,\gd v, r_4, \ga l_1]
 =\Bigl[B(x),\frac{B(y)}{(b-1)y}, \frac{abx+1}{abA(x)},\frac{(b-1)(y-x)}{xB(y)},\ga l_1(y)\Bigr]$}}} $$
The same proof above works if we notice now that
$$
 \pa_1^0[B,B/f]=\Bigl[\frac{B(y)}{(b-1)y}, \frac{1}{b},b-1,\ga l_1(y)\Bigr]
$$
is admissible.

$$\text{\boxed{\text{$[B,B/f]_2:=[gh, g,  v, r_4, \ga l_1]
 =\Bigl[B(x),\frac{B(y)}{(b-1)y}, \frac{abx+1}{aA(x)},\frac{(b-1)(y-x)}{xB(y)},\ga l_1(y)\Bigr]$}}} $$
We have
\begin{alignat*}{3}
 \pa_1^0[B,B/f]_2& \subset \{t_3=1\},\quad&
   \pa_2^0[B,B/f]_2&\subset \{t_5=1\}, \\
\pa_4^\infty[B,B/f]_2&\subset\{t_5=1\}, \quad &
 \pa_5^\infty[B,B/f]_2&\subset\{t_2=1\},
\end{alignat*}
and
\begin{align*}
\pa_1^\infty[B,B/f]_2&=\Bigl[\frac{B(y)}{(b-1)y},b,\frac{1-b}{B(y)},\ga l_1(y)\Bigr],\\
 \pa_2^\infty[B,B/f]_2&=\Bigl[B(x),\frac{abx+1}{aA(x)},1-b,\ga \Bigr],\\
 \pa_3^0[B,B/f]_2&=\Bigl[\frac{\mu}{b},\frac{B(y)}{(b-1)y},
       \frac{(1-b)(aby+1)}{B(y)},\ga l_1(y) \Bigr],\\
 \pa_3^\infty[B,B/f]_2&=\Bigl[-\mu,\frac{B(y)}{(b-1)y},
       \frac{a(b-1)A(y)}{(a-1)B(y)},\ga l_1(y) \Bigr],\\
 \pa_4^0[B,B/f]_2&=\Bigl[B(y),\frac{B(y)}{(b-1)y},
   \frac{aby+1}{aA(y)},\ga l_1(y) \Bigr],\\
  \pa_5^0[B,B/f]_2&=\Bigl[B(x),\frac{B(c)}{(b-1)c},\frac{abx+1}{abA(x)},
   \frac{(b-1)(c-x)}{xB(c)}\Bigr].
\end{align*}
All these cycles are clearly admissible by our choice of $\ga$.
The same proof shows that $[gh, g,  \gd, r_4, \ga l_1]$ is
admissible.

$$\text{\boxed{\text{$[B,B/f]_2':=[gh, g,  v, r_4/w_4, \ga l_1]
 =\Bigl[B(x),\frac{B(y)}{(b-1)y}, \frac{abx+1}{aA(x)},\frac{(b-1)B(x)(y-1)}{xB(y)},\ga l_1(y)\Bigr]$}}} $$
This is a sum of two admissible cycles:
$$\Bigl[B(x),\frac{B(y)}{(b-1)y},
\frac{abx+1}{aA(x)},\frac{(y-1)}{B(y)},\ga l_1(y)\Bigr]
 +\Bigl[B(x),\frac{B(y)}{(b-1)y},
 \frac{abx+1}{aA(x)},\frac{(b-1)B(x)}{x},\ga l_1(y)\Bigr].$$

$$\text{\boxed{\text{$[B,B/f]_3:=[gh, g,v, w_4, \ga l_1]
 =\Bigl[B(x),\frac{B(y)}{(b-1)y}, \frac{abx+1}{aA(x)},\frac{y-x}{B(x)(y-1)},\ga l_1(y)\Bigr]$}}} $$
We have
$$ \pa_1^0[B,B/f]_3 \subset \{t_3=1\},\quad
  \pa_2^0[B,B/f]_3\subset \{t_5=1\}, \quad
 \pa_5^\infty[B,B/f]_3\subset\{t_2=1\},$$
and
\begin{align*}
\pa_1^\infty[B,B/f]_3&=\Bigl[\frac{B(y)}{(b-1)y},b,\frac{y}{(1-b)(y-1)},\ga l_1(y)\Bigr],\\
 \pa_2^\infty[B,B/f]_3&=\Bigl[B(x),\frac{abx+1}{abA(x)},\frac{x}{B(x)},\ga \Bigr],\\
 \pa_3^0[B,B/f]_3&=\Bigl[\frac{\mu}{b},\frac{B(y)}{(b-1)y},
       \frac{a(1-b)(aby+1)}{\mu(y-1)},\ga l_1(y) \Bigr],\\
 \pa_3^\infty[B,B/f]_3&=\Bigl[-\mu,\frac{B(y)}{(b-1)y},
       \frac{A(y)}{\mu(1-y)},\ga l_1(y) \Bigr],\\
 \pa_4^0[B,B/f]_3&=\Bigl[B(y),\frac{B(y)}{(b-1)y},
   \frac{aby+1}{aA(y)},\ga l_1(y) \Bigr],\\
\pa_4^\infty[B,B/f]_3&=
 =\Bigl[B(x),\frac{1}{b-1}, \frac{abx+1}{aA(x)},\frac{\ga
 (c-1)}{c}\Bigr]\\
  \pa_5^0[B,B/f]_3&=\Bigl[B(x),\frac{B(c)}{(b-1)c},\frac{abx+1}{aA(x)},
   \frac{c-x}{B(x)(c-1)}\Bigr].
\end{align*}
All these cycles are clearly admissible.

$$\text{\boxed{\text{$[B,B/f]_4:=[gh, g,v, w_4, l_1]
 =\Bigl[B(x),\frac{B(y)}{(b-1)y}, \frac{abx+1}{aA(x)},\frac{y-x}{B(x)(y-1)},l_1(y)  \Bigr]$}}} $$
Note that in the above proof for $[B,B/f]_3$ the choice of $\ga$
is not essential because whenever $B(x)=0$ we have
$(abx+1)/aA(x)=1$. The same reason shows that $[gh,g,v,w_4, \ga]$
is admissible.

$$\text{\boxed{\text{$[B,B/f]_4':=[gh, g,v, w_4/p_4, l_1]
 =\Bigl[B(x),\frac{B(y)}{(b-1)y}, \frac{abx+1}{aA(x)},\frac{A(y)}{1-y},l_1(y)\Bigr]$}}} $$
This is a product of two admissible cycles.

$$\text{\boxed{\text{$[B,B/f]_5:=[gh, g,v, p_4, l_1]
 =\Bigl[B(x),\frac{B(y)}{(b-1)y}, \frac{abx+1}{aA(x)},\frac{\mu(x-y)}{A(y)B(x)},l_1(y)\Bigr]$}}} $$
The proof for $[B,B/f]_4$ can be adapted here without any change.

\item $[gh, h,\cdots]$ are admissible.

$$\text{\boxed{\text{$[B,f]:=[gh, h,\gd v, s_4, \ga l_1]
 =\Bigl[B(x),(b-1)y, \frac{abx+1}{abA(x)},\frac{(b-1)(y-x)}{B(y)},\ga l_1(y)\Bigr]$}}} $$
We have
\begin{alignat*}{3}
 \pa_1^0[B,f]&\subset \{t_4=1\},\quad&
 \pa_1^\infty[B,f]&\subset \{t_3=1\},\quad
  \pa_2^\infty[B,f]\subset \{t_4=1\},\\
 \pa_4^\infty[B,f]&\subset \{t_5=1\},\quad&
 \pa_5^\infty[B,f]&\subset\{t_4=1\},
\end{alignat*}
and
\begin{align*}
 \pa_2^0[B,f]&=\Bigl[B(x),\frac{abx+1}{abA(x)},(1-b)x,\ga \Bigr],\\
 \pa_3^0[B,B/f]&=\Bigl[\frac{\mu}{b},(b-1)y,
       \frac{(b-1)(aby+1)}{abB(y)},\ga l_1(y) \Bigr],\\
 \pa_3^\infty[B,B/f]&=\Bigl[-\mu,(b-1)y,
       \frac{(b-1)A(y)}{B(y)},\ga l_1(y) \Bigr],\\
 \pa_4^0[B,f]&=\Bigl[B(y),(b-1)y,
    \frac{aby+1}{abA(y)},\ga l_1(y) \Bigr],\\
 \pa_5^0[B,f]&=\Bigl[B(x),(1-b)c,\frac{abx+1}{abA(x)},\frac{(b-1)(c-x)}{B(c)}
    \Bigr].
\end{align*}
All these cycles are clearly admissible. Note that we never use
the property of $\ga$, (namely $t_5=1$), so the same proof shows
that $[gh, h,\gd v, q_4, l_1]$ is admissible.

$$\text{\boxed{\text{$[B,f]':=[gh, h,\gd v, s_4/q_4, l_1]
 =\Bigl[B(x),(b-1)y, \frac{abx+1}{abA(x)},\frac{(b-1)A(y)}{B(y)}, l_1(y)\Bigr]$}}} $$
This is a product of two admissible cycles.

$$\text{\boxed{\text{$[B,f]_1:=[gh, h,\gd v, q_4, l_1]
 =\Bigl[B(x),(b-1)y, \frac{abx+1}{abA(x)},\frac{y-x}{A(y)}, l_1(y)\Bigr]$}}} $$
The same proof for $[B,f]$ works.

\item $[h,gh, \cdots]$ are admissible.

$$\text{\boxed{\text{$[f,B]:=[h,gh,\gd v, q_4,\ga l_1]
 =\Bigl[(b-1)x,B(y), \frac{abx+1}{abA(x)},\frac{y-x}{A(y)},\ga l_1(y)\Bigr]$}}} $$
We have
\begin{alignat*}{3}
 \pa_1^\infty[f,B]&\subset \{t_3=1\} \Bigr],\quad &
 \pa_2^0[f,B]&\subset \{t_5=1\},\quad
  \pa_2^\infty[f,B] \subset \{t_4=1\}, \\
 \pa_3^\infty[f,B]&\subset \{t_4=1\},\quad &
 \pa_5^\infty[f,B]&\subset\{t_4=1\},
\end{alignat*}
and
\begin{align*}
 \pa_1^0[f,B]&=\Bigl[B(y),\frac{1}{b},
    \frac{y}{A(y)},\ga l_1(y) \Bigr],\\
 \pa_3^0[f,B]&=\Bigl[\frac{1-b}{ab},B(y),
       \frac{aby+1}{abA(y)},\ga l_1(y) \Bigr],\\
 \pa_4^0[f,B]&=\Bigl[(b-1)y,B(y),
    \frac{aby+1}{abA(y)},\ga l_1(y) \Bigr],\\
 \pa_4^\infty[f,B]&=\Bigl[(b-1)x,-\mu,
    \frac{abx+1}{abA(x)}, \frac{(b-1)(ac-a+1)}{a(bc-c+1)}\Bigr],\\
 \pa_5^0[f,B]&=\Bigl[(b-1)x,B(c),\frac{abx+1}{abA(x)},
    \frac{c-x}{A(c)}\Bigr].
\end{align*}
All these cycles are clearly admissible by our choice of $\gd$.

$$\text{\boxed{\text{$[f,B]_1:=[h,gh,\gd v, q_4,l_1]
 =\Bigl[(b-1)x,B(y), \frac{abx+1}{abA(x)},\frac{y-x}{A(y)},l_1(y)\Bigr]$}}} $$

Note that we use the property of $\ga$ (namely $t_5=1$) only for
$\pa_2^0[f,B]$ so the same proof applies because
$$\pa_2^0[f,B]_1=\Bigl[(b-1)x,\frac{abx+1}{abA(x)},
 \frac{B(x)}{-\mu},\frac{1}{\ga}\Bigr]$$
is clearly admissible.

\item $[g, h,\cdots]$ are admissible.

$$\text{\boxed{\text{$[B/f,f]:=[g, h,\gd v, q_4, l_1]
 =\Bigl[\frac{B(x)}{(b-1)x},(b-1)y, \frac{abx+1}{abA(x)},\frac{y-x}{A(y)}, l_1(y)\Bigr]$}}} $$
We have
\begin{alignat*}{3}
 \pa_2^0[B/f,f]&\subset \{t_5=1\},\quad&
  \pa_2^\infty[B/f,f]&\subset \{t_4=1\},  \\
 \pa_3^\infty[B/f,f]&\subset \{t_4=1\},\quad&
 \pa_5^\infty[B/f,f]&\subset\{t_4=1\},
\end{alignat*}
and
\begin{align*}
 \pa_1^0[B/f,f]&=\Bigl[(b-1)y,\frac{1}{b},
    \frac{B(y)}{(b-1)A(y)}, l_1(y) \Bigr],\\
 \pa_1^\infty[B/f,f]&=\Bigl[(b-1)y,\frac{1}{b(1-a)},
    \frac{y}{A(y)},l_1(y)  \Bigr],\\
 \pa_3^0[B/f,f]&=\Bigl[\frac{a\mu}{1-b},(b-1)y,
       \frac{aby+1}{abA(y)}, l_1(y) \Bigr],\\
 \pa_4^0[B/f,f]&=\Bigl[\frac{B(y)}{(b-1)y},(b-1)y,
    \frac{aby+1}{abA(y)},l_1(y) \Bigr],\\
 \pa_4^\infty[B/f,f]&=\Bigl[\frac{B(x)}{(b-1)x},\frac{(a-1)(b-1)}{a},
    \frac{abx+1}{abA(x)}, \frac{ac-a+1}{ac}\Bigr],\\
 \pa_5^0[B/f,f]&=\Bigl[\frac{B(x)}{(b-1)x},(1-b)c,\frac{abx+1}{abA(x)},
    \frac{c-x}{A(c)}\Bigr].
\end{align*}
All these cycles are clearly admissible. Note that we never use
the property of $\gd$, (namely $t_3=1$), so the same proof shows
that $[g, h, v, q_4, l_1]$ and $[g, h, \gd, q_4, l_1]$ are
admissible.

\item $[h,g, \cdots]$ are admissible.

$$\text{\boxed{\text{$[f,B/f]:=[h,g, \gd v, q_4, l_1]
 =\Bigl[(b-1)x,\frac{B(y)}{(b-1)y},\frac{abx+1}{abA(x)},\frac{y-x}{A(y)}, l_1(y)\Bigr]$}}} $$
We have
\begin{alignat*}{3}
 \pa_1^\infty[f,B/f]&\subset \{t_3=1\},\quad&
  \pa_2^\infty[f,B/f]&\subset \{t_5=1\},  \\
 \pa_3^\infty[f,B/f]&\subset \{t_4=1\},\quad&
 \pa_5^\infty[f,B/f]&\subset\{t_2=1\},
\end{alignat*}
and
\begin{align*}
 \pa_1^0[f,B/f]&=\Bigl[\frac{B(y)}{(b-1)y},\frac{1}{b(1-a)},
    \frac{y}{A(y)},l_1(y)  \Bigr],\\
 \pa_2^0[f,B/f]& =\Bigl[(b-1)x,\frac{abx+1}{abA(x)},\frac{B(x)}{-\mu},
    \frac{1}{\ga}\Bigr]\\
 \pa_3^0[f,B/f]&=\Bigl[\frac{1-b}{ab},\frac{B(y)}{(b-1)y},
       \frac{aby+1}{abA(y)}, l_1(y) \Bigr],\\
 \pa_4^0[f,B/f]&=\Bigl[(b-1)y,\frac{B(y)}{(b-1)y},
    \frac{aby+1}{abA(y)},l_1(y) \Bigr],\\
 \pa_4^\infty[f,B/f]&=\Bigl[(b-1)x,\frac{ab-b+1}{(a-1)(b-1)},
    \frac{abx+1}{abA(x)}, \frac{ac-a+1}{ac}\Bigr],\\
 \pa_5^0[f,B/f]&=\Bigl[(b-1)x,\frac{B(c)}{(b-1)c},\frac{abx+1}{abA(x)},
    \frac{c-x}{A(c)}\Bigr].
\end{align*}
All these cycles are clearly admissible.

$$\text{\boxed{\text{$[f,B/f]':=[h,g, \gd v, p_4/q_4, l_1]
 =\Bigl[(b-1)x,\frac{B(y)}{(b-1)y},\frac{abx+1}{abA(x)},
 \frac{-\mu}{B(x)}, l_1(y)\Bigr]$}}} $$
In the above proof the only place we use $q_4=1$ (namely $t_4=1$)
is for $\pa_3^\infty[f,B/f]$. However, it's still true in
$[f,B/f]_1'$ that $t_4=1$ if $A(x)=0$. Now the only things that
need checking are
$$\pa_4^0[f,B/f]'\subset \{t_3=1\},\quad
 \pa_4^\infty[f,B/f]'=\Bigl[-1,\frac{B(y)}{(b-1)y},\frac{1}{b},l_1(y)\Bigr].$$

$$\text{\boxed{\text{$[f,B/f]_1:=[h,g, \gd v, p_4, l_1]
 =\Bigl[(b-1)x,\frac{B(y)}{(b-1)y},\frac{abx+1}{abA(x)},\frac{\mu(x-y)}{A(y)B(x)}, l_1(y)\Bigr]$}}} $$
Similar to $[f,B/f]'$ the only things that need checking are
 $$\pa_4^0[f,B/f]_1=\pa_4^0[f,B/f],\quad
\pa_4^\infty[f,B/f]_1=\pa_4^\infty[f,B/f]+\pa_4^\infty[f,B/f]_1'.$$

$$\text{\boxed{\text{$[f,B/f]_2:=[h,g, v, p_4, l_1]
 =\Bigl[(b-1)x,\frac{B(y)}{(b-1)y},\frac{abx+1}{aA(x)},\frac{\mu(x-y)}{A(y)B(x)}, l_1(y)\Bigr]$}}} $$

In the above proofs for $[f,B/f]$ and $[f,B/f]_1$ the only place
we use the property of $\gd$ (namely $t_3=1$) for
$\pa_1^\infty[f,B/f]$. But
$$\pa_1^\infty[f,B/f]_2=\Bigl[\frac{B(y)}{(b-1)y},
       \frac{1}{b},\frac{\mu}{(b-1)A(y)}, l_1(y) \Bigr]$$
which is admissible. This shows that $[f,B/f]_2$ is admissible.

\item $[g,g, \cdots]$ are admissible.

$$\text{\boxed{\text{$[B/f,B/f]:=[g,g, v, q_4, l_1]
 =\Bigl[\frac{B(x)}{(b-1)x},\frac{B(y)}{(b-1)y},\frac{abx+1}{aA(x)},\frac{y-x}{A(y)}, l_1(y)\Bigr]$}}} $$
We have
\begin{alignat*}{3}
 \pa_1^0[B/f,B/f]&\subset \{t_3=1\},\quad&
  \pa_2^\infty[B/f,B/f]&\subset \{t_5=1\},  \\
 \pa_3^\infty[B/f,B/f]&\subset \{t_4=1\},\quad&
 \pa_5^\infty[B/f,B/f]&\subset\{t_2=1\},
\end{alignat*}
and
\begin{align*}
 \pa_1^\infty[B/f,B/f]&=\Bigl[\frac{B(y)}{(b-1)y},\frac{1}{1-a},
    \frac{y}{A(y)},l_1(y)\Bigr],\\
 \pa_2^0[B/f,B/f]& =\Bigl[\frac{B(x)}{(b-1)x},\frac{abx+1}{aA(x)},\frac{B(x)}{-\mu},
    \frac{1}{\ga}\Bigr]\\
 \pa_3^0[B/f,B/f]&=\Bigl[\frac{a\mu}{b-1},\frac{B(y)}{(b-1)y},
       \frac{aby+1}{abA(y)}, l_1(y) \Bigr],\\
 \pa_4^0[B/f,B/f]&=\Bigl[\frac{B(y)}{(b-1)y},\frac{B(y)}{(b-1)y},
    \frac{aby+1}{aA(y)},l_1(y) \Bigr],\\
 \pa_4^\infty[B/f,B/f]&=\Bigl[\frac{B(x)}{(b-1)x},\frac{ab-b+1}{(a-1)(b-1)},
    \frac{abx+1}{aA(x)}, \frac{ac-a+1}{ac}\Bigr],\\
 \pa_5^0[B/f,B/f]&=\Bigl[\frac{B(x)}{(b-1)x},\frac{B(c)}{(b-1)c},\frac{abx+1}{aA(x)},
    \frac{c-x}{A(c)}\Bigr].
\end{align*}
All these cycles are clearly admissible.

$$\text{\boxed{\text{$[B/f,B/f]_1:=[g,g, v, p_4, l_1]
 =\Bigl[\frac{B(x)}{(b-1)x},\frac{B(y)}{(b-1)y},\frac{abx+1}{aA(x)},
 \frac{\mu(x-y)}{A(y)B(x)}, l_1(y)\Bigr]$}}} $$
In the above proof the only place we use $q_4=1$ (namely $t_4=1$)
is for $\pa_3^\infty[B/f,B/f]$. However, it's still true in
$[B/f,B/f]_1$ that $p_4=1$ if $A(x)=0$. Now the only things that
need checking are
\begin{align*}
\pa_4^0[B/f,B/f]_1=&\pa_4^0[B/f,B/f],\\
\pa_4^\infty[B/f,B/f]_1=&\pa_4^\infty[B/f,B/f]\quad \text{(if
 $B(x)=0$   then $t_3=1$}).
\end{align*}

$$\text{\boxed{\text{$[B/f,B/f]_1':=[g,g, v, p_4/q_4, l_1]
 =\Bigl[\frac{B(x)}{(b-1)x},\frac{B(y)}{(b-1)y},\frac{abx+1}{aA(x)},
 \frac{-\mu}{B(x)}, l_1(y)\Bigr]$}}} $$
In the above proof the only place we use $q_4=1$ (namely $t_4=1$)
is for $\pa_3^\infty[B/f,B/f]$. However, it's still true in
$[B/f,B/f]_3$ that $t_4=1$ if $A(x)=0$. Now the only things that
need checking are
$$\pa_4^0[B/f,B/f]_1'\subset \{t_1=1\},\quad
 \pa_4^\infty[B/f,B/f]_1'\subset \{t_3=1\}.$$

\item $[h,h,\gd v, q_4, l_1]$ is admissible.
$$\text{\boxed{\text{$[f,f]:=[h,h, \gd v, q_4, l_1]
 =\Bigl[(b-1)x,(b-1)y,\frac{abx+1}{abA(x)},\frac{y-x}{A(y)}, l_1(y)\Bigr]$}}} $$
We have
\begin{alignat*}{3}
 \pa_1^\infty[f,f]&\subset \{t_3=1\},\quad&
 \pa_2^0[f,f]&\subset \{t_5=1\},\quad
  \pa_2^\infty[f,f]\subset \{t_4=1\},  \\
 \pa_3^\infty[f,f]&\subset \{t_4=1\},\quad&
 \pa_5^\infty[f,f]&\subset\{t_4=1\},
\end{alignat*}
and
\begin{align*}
 \pa_1^0[f,f]&=\Bigl[(b-1)y,\frac{1}{b(1-a)},
    \frac{y}{A(y)},l_1(y)  \Bigr],\\
 \pa_3^0[f,f]&=\Bigl[\frac{1-b}{ab},(b-1)y,
       \frac{aby+1}{abA(y)}, l_1(y) \Bigr],\\
 \pa_4^0[f,f]&=\Bigl[(b-1)y,(b-1)y,
    \frac{aby+1}{abA(y)},l_1(y) \Bigr],\\
 \pa_4^\infty[f,f]&=\Bigl[(b-1)x,\frac{(a-1)(b-1)}{a},
    \frac{abx+1}{abA(x)}, \frac{ac-a+1}{ac}\Bigr],\\
 \pa_5^0[f,f]&=\Bigl[(b-1)x,(1-b)c,\frac{abx+1}{abA(x)},
    \frac{c-x}{A(c)}\Bigr].
\end{align*}
All these cycles are clearly admissible. Note that we didn't use
$t_1=1$ or $t_2=1$ in the above so the same argument implies that
$$[b-1,h, \gd v, q_4, l_1],\ [h,b-1, \gd v, q_4, l_1],
    \ [b-1,b-1, \gd v, q_4,l_1]$$
are all admissible. So we get
\begin{equation}\label{h9}
[h,h, \gd v,q_4,l_1] =\Bigl[x,y,\frac{abx+1}{abA(x)},
 \frac{y-x}{A(y)}, l_1(y)\Bigr]
\end{equation}

\end{itemize}

\medskip
\noindent {\bf Step} (8). Computation of $Y_3+Y_4$.
\medskip

\medskip
\noindent{\bf Claim 8.1.}  Under non-degeneracy assumption
$$Y'_{31}:=\Bigl[\frac{(1-b)x}{B(x)},(1-b)y,\frac{abA(x)}{abx+1},
 \frac{ab(y-x)}{aby+1},l_2(y)\Bigr]=-Y_3'.$$

$$\text{\boxed{\text{$Y'_{31}=\Bigl[\frac{(1-b)x}{B(x)}, (1-b)y,\frac{aA(x)}{abx+1},
 \frac{ab(y-x)}{aby+1},l_2(y)\Bigr]$}}} $$
We have
\begin{alignat*}{3}
 \pa_1^\infty Y'_{31}&\subset \{t_3=1\},\quad&
 \pa_2^0 Y'_{31}&\subset \{t_5=1\},\quad
  \pa_2^\infty Y'_{31}\subset \{t_4=1\},  \\
 \pa_3^\infty Y'_{31}&\subset \{t_4=1\},\quad&
 \pa_5^\infty Y'_{31}&\subset\{t_4=1\},
\end{alignat*}
and
\begin{align*}
 \pa_1^0 Y'_{31}&=\Bigl[(1-b)y,1-a,\frac{aby}{aby+1},l_2(y)\Bigr],\\
 \pa_3^0 Y'_{31}&=\Bigl[\frac{(a-1)(1-b)}{ab-b+1},(1-b)y,
       \frac{abA(y)}{aby+1}, l_2(y) \Bigr],\\
 \pa_4^0 Y'_{31}&=\Bigl[\frac{(1-b)y}{B(y)}, (1-b)y,\frac{aA(y)}{aby+1},l_2(y)\Bigr]],\\
 \pa_4^\infty Y'_{31}&=\Bigl[\frac{(1-b)x}{B(x)},\frac{b-1}{ab},
    \frac{aA(x)}{abx+1}, \frac{ac-a}{ac-a+1}\Bigr],\\
 \pa_5^0 Y'_{31}&=\Bigl[\frac{(1-b)x}{B(x)}, (1-b)y_2,\frac{aA(x)}{abx+1},
 \frac{ab(y_2-x)}{aby_2+1}\Bigr].
\end{align*}
All these cycles are clearly admissible. For the last one, we need
\eqref{aby2}.

$$\text{\boxed{\text{$Y''_{31}=\Bigl[\frac{(1-b)x}{B(x)}, (1-b)y,b,
 \frac{ab(y-x)}{aby+1},l_2(y)\Bigr]$}}} $$
In the above proof for $Y'_{31}$ there is only one place where we
used $t_3=1$, namely for $\pa_1^\infty$. But
$$ \pa_3^\infty Y''_{31}=\Bigl[(1-b)y,b,
 \frac{abB(y)}{(b-1)(aby+1)},l_2(y)\Bigr]$$
 which is admissible because if $B(y)=0$ then $(1-b)y=1$.

$$\text{\boxed{\text{$Y'_{32}=\Bigl[\frac{(1-b)x}{B(x)}, (1-b)y,\frac{abA(x)}{abx+1},
\frac{aby+1}{abA(y)},l_2(y)\Bigr]$}}}
$$
This is  a product of two admissible cycles. From the above we get
$$Y_{31}'=Y_{31}'+Y_{32}'=Y_{31}'+Y_{32}'+Y_{31}''=-Y_3'.$$
Claim 8.1 is proved.

\medskip
\noindent{\bf Claim 8.2.} Under non-degeneracy assumption
$$Y'_{41}:=\Bigl[(1-b)x,\frac{(1-b)y}{B(y)},\frac{abA(x)}{abx+1},
 \frac{(ab-b+1)(y-x)}{(aby+1)B(x)},l_2(y)\Bigr]=-Y_4'.$$

$$\text{\boxed{\text{$Y'_{41}=\Bigl[(1-b)x,\frac{B(y)}{(1-b)y},\frac{aA(x)}{abx+1},
 \frac{(ab-b+1)(y-x)}{(aby+1)B(x)},l_2(y)\Bigr]$}}} $$
We have
 $$ \pa_2^0  Y'_{41}\subset \{t_4=1\},\quad
  \pa_2^\infty  Y'_{41}\subset \{t_5=1\},\quad
 \pa_5^\infty  Y'_{41}\subset\{t_4=1\},   \quad
 \pa_3^\infty  Y'_{41}\subset \{t_4=1\},
$$and
\begin{align*}
 \pa_1^\infty  Y'_{41}&=\Bigl[\frac{B(y)}{(1-b)y},1-a,\frac{(ab-b+1)y}{aby+1},l_2(y)\Bigr],\\
 \pa_1^0  Y'_{41}&=\Bigl[\frac{B(y)}{(1-b)y},\frac{1}{b},\frac{ab-b+1}{(1-b)(aby+1)},l_2(y)\Bigr],\\
 \pa_3^0  Y'_{41}&=\Bigl[\frac{(a-1)(1-b)}{a},\frac{B(y)}{(1-b)y},
       \frac{aA(y)}{aby+1}, l_2(y) \Bigr],\\
 \pa_4^0  Y'_{41}&=\Bigl[(1-b)y,\frac{B(y)}{(1-b)y}, \frac{aA(y)}{aby+1},l_2(y)\Bigr]],\\
 \pa_4^\infty Y'_{41}&=\Bigl[(1-b)x,\frac{ab-b+1}{b-1},\frac{aA(x)}{abx+1}, \frac{ac-a}{ac-a+1} \Bigr],\\
 \pa_5^0  Y'_{41}&=\Bigl[(1-b)x,\frac{B(y_2)}{(1-b)y_2}, \frac{aA(x)}{abx+1},
 \frac{(ab-b+1)(y_2-x)}{(aby_2+1)B(x)}\Bigr].
\end{align*}
All these cycles are clearly admissible. For the last one, we need
\eqref{aby2} and \eqref{By2}.

$$\text{\boxed{\text{$Y''_{41}=\Bigl[(1-b)x,\frac{B(y)}{(1-b)y},b,
 \frac{(ab-b+1)(y-x)}{(aby+1)B(x)},l_2(y)\Bigr]$}}} $$
In the above proof for $Y'_{41}$ we didn't use $t_3=1$ so the same
proof is still valid.

$$\text{\boxed{\text{$Y'_{42}=\Bigl[(1-b)x,\frac{B(y)}{(1-b)y},\frac{abA(x)}{abx+1},
 \frac{aby+1}{aA(y)},l_2(y)\Bigr]$}}} $$
 This is  a product of two admissible cycles. From the above we get
$$Y_{41}'=Y_{41}'+Y_{42}'=Y_{41}'+Y_{42}'+Y_{41}''=-Y_4'.$$
Claim 8.2 is proved.

\medskip
\noindent{\bf Claim 8.3.}  Under non-degeneracy assumption
$$Y_3'+Y_4'=Y_3+Y_4.$$

First it is not hard to show all the following cycles are
admissible and negligible:
\begin{align*}
Y_{31}=&\Bigl[-1, (b-1)y,\frac{abx+1}{aA(x)},
    \frac{y-x}{A(y)},l_2(y)\Bigr],\\
Y_{41}=&\Bigl[(b-1)x,-1,\frac{abx+1}{aA(x)},
    \frac{\mu(x-y)}{A(y)B(x)},l_2(y)\Bigr],\\
Y_{32}=&\Bigl[\frac{B(x)}{(1-b)x}, -1,\frac{abx+1}{aA(x)},
    \frac{y-x}{A(y)},l_2(y)\Bigr],\\
Y_{42}=&\Bigl[-1,\frac{B(y)}{(1-b)y},\frac{abx+1}{aA(x)},
    \frac{\mu(x-y)}{A(y)B(x)},l_2(y)\Bigr].
\end{align*}
Then by Lemma 3.2(ii) we have
\begin{align*}
Y_3'+Y_4'=&\Bigl[\frac{B(x)}{(1-b)x},(b-1)y,\frac{abx+1}{aA(x)},
    \frac{\mu(x-y)}{A(y)B(x)},l_2(y)\Bigr]+Y_{32}\\
    +&\Bigl[(b-1)x,\frac{B(y)}{(1-b)y},\frac{abx+1}{aA(x)},
    \frac{\mu(x-y)}{A(y)B(x)},l_2(y)\Bigr]+Y_{42}\\
    =&Y_3+Y_4+Y_{31}+Y_{41}\\
    =&Y_3+Y_4.
\end{align*}
Claim 8.3 is proved.

\medskip
\noindent {\bf Step} (9). Final decomposition of $\{k(c)\}$ into
$T_i(F)$'s.

$$\text{\boxed{\text{$T_1(A)=\Bigl[A(x),A(y),
    \frac{x-1}{x}, \frac{y-x}{A(y)},\eps_1(A) l_1(y)\Bigr],\quad \eps_1(A)=\frac{ac}{ac-a+1}$}}} $$
We have
\begin{alignat*}{3}
 \pa_1^0 T_1(A)&\subset \{t_4=1\},\quad&
 \pa_1^\infty T_1(A)&\subset \{t_3=1\},\quad
 \pa_2^0 T_1(A) \subset \{t_5=1\},\\
 \pa_2^\infty T_1(A)&\subset \{t_4=1\},\quad&
 \pa_4^\infty T_1(A)&\subset \{t_5=1\},   \quad
 \pa_5^\infty T_1(A)\subset\{t_4=1\},
\end{alignat*}
and
\begin{align*}
 \pa_3^0  T_1(A)&=\Bigl[\frac{1}{a},A(y),
       \frac{y-1}{A(y)},\eps_1(A) l_1(y) \Bigr],\\
 \pa_3^\infty T_1(A)&=\Bigl[\frac{1-a}{a},A(y),
       \frac{y}{A(y)},\eps_1(A) l_1(y) \Bigr],\\
 \pa_4^0  T_1(A)&=\Bigl[A(y),A(y),
    \frac{y-1}{y}, \eps_1(A) l_1(y) \Bigr],\\
 \pa_5^0  T_1(A)&=\Bigl[A(x),A(c),
    \frac{x-1}{x}, \frac{c-x}{A(c)} \Bigr].
\end{align*}
All these cycles are clearly admissible.

$$\text{\boxed{\text{$T_2(A)=\Bigl[A(x),A(y),
   \frac{x-1}{x}, \frac{y-x}{A(y)},\eps_2(A) l_2(y)\Bigr],\quad \eps_2(A)=\frac{ac-a+1}{ac} $}}} $$
The above proof mostly is still valid because $\eps_2(A)
l_2((a-1)/a)=1$ except that
\begin{align*}
 \pa_5^0  T_2(A)&=\Bigl[A(x),A(y_2),
    \frac{x-1}{x}, \frac{y_2-x}{A(y_2)} \Bigr],\\
\pa_5^\infty T_2(A)&=\Bigl[A(x),\frac{\mu}{b-1},
    \frac{x-1}{x},  \frac{B(x)}{-\mu} \Bigr]
\end{align*}
which are both admissible.

$$\text{\boxed{\text{$T_3(A)=\Bigl[A(x),A(y), \frac{abx+1}{abA(x)},
    \frac{y-x}{A(y)}, \eps_1(A)l_1(y)\Bigr],\quad \eps_1(A)=\frac{ac}{ac-a+1} $}}} $$

We have
\begin{alignat*}{3}
 \pa_1^0 T_3(A)&\subset \{t_4=1\},\quad&
 \pa_1^\infty T_3(A)&\subset \{t_3=1\},\quad
  \pa_2^0 T_3(A)  \subset \{t_5=1\},\quad
 \pa_2^\infty T_3(A)\subset \{t_4=1\},  \\
 \pa_3^\infty T_3(A)&\subset \{t_4=1\},\quad&
 \pa_4^\infty T_3(A)&\subset \{t_5=1\},   \quad
 \pa_5^\infty T_3(A)\subset\{t_4=1\},
\end{alignat*}
and
\begin{align*}
 \pa_3^0  T_3(A)&=\Bigl[\frac{\mu}{b},A(y),
       \frac{aby+1}{abA(y)},\eps_1 l_1(y) \Bigr],\\
 \pa_4^0  T_3(A)&=\Bigl[A(y),A(y),
    \frac{aby+1}{abA(y)},\eps_1 l_1(y) \Bigr],\\
 \pa_5^0  T_3(A)&=\Bigl[A(x),A(c),
    \frac{abx+1}{abA(x)}, \frac{c-x}{A(c)} \Bigr].
\end{align*}
All these cycles are clearly admissible.

$$\text{\boxed{\text{$T_4(A)=\Bigl[A(x),A(y), \frac{abx+1}{abA(x)},
    \frac{y-x}{A(y)}, \eps_2(A)l_1(y)\Bigr],\quad \eps_2(A)=\frac{ac-a+1}{ac} $}}} $$
The above proof mostly is still valid because $\eps_2
l_2((a-1)/a)=1$ except that
\begin{align*}
 \pa_5^0  T_4(A)&=\Bigl[A(x),A(y_2),
    \frac{abx+1}{abA(x)}, \frac{y_2-x}{A(y_2)} \Bigr],\\
\pa_5^\infty T_4(A)&=\Bigl[A(x),\frac{\mu}{b-1},
   \frac{abx+1}{abA(x)}, \frac{B(x)}{-\mu} \Bigr]
\end{align*}
which are both admissible.

Next we need to prove
\begin{align*}
Z_3\Bigl(\frac{A}{f},\frac{A}{f}\Bigr)
=&\Bigl[\frac{A(x)}{x},\frac{A(y)}{y},
    \frac{x-1}{x}, \frac{y-x}{yB(x)},l(y)\Bigr]
=\Bigl[\frac{A(x)}{x},\frac{A(y)}{y},
    \frac{x-1}{x}, 1-\frac{x}{y},l(y)\Bigr]\\
=&\Bigl[\frac{A(x)}{x},\frac{A(y)}{y},
    \frac{(1-a)(x-1)}{x}, \frac{y-x}{yB(x)},l(y)\Bigr].
\end{align*}

$$\text{\boxed{\text{$Z_3'(A)=\Bigl[\frac{A(x)}{x},\frac{A(y)}{y},
    \frac{x-1}{x}, \frac{y-x}{y},l(y)\Bigr]$}}} $$
The only non-trivial boundaries are
\begin{align*}
  \pa_3^0 Z_3'(A) &=\Bigl[\frac{1-a}{a},\frac{A(y)}{y},
       \frac{y-1}{y},l(y) \Bigr],\\
 \pa_4^0  Z_3'(A) &=[\frac{A(y)}{y},\frac{A(y)}{y},
    \frac{y-1}{y}, l(y)\Bigr],\\
 \pa_5^0  Z_3'(A)&=\Bigl[\frac{A(x)}{x},\frac{A(c)}{c},
    \frac{x-1}{x}, \frac{c-x}{c}\Bigr]
    +\Bigl[\frac{A(x)}{x},\frac{A(y_2)}{y_2},
    \frac{x-1}{x}, \frac{y_2-x}{y_2}\Bigr]
\end{align*}
which are all admissible.

$$\text{\boxed{\text{$Z_3''(A)=\Bigl[\frac{A(x)}{x},\frac{A(y)}{y},
    \frac{x-1}{x}, \frac{1}{B(x)},l(y)\Bigr]$}}} $$
It's admissible because $B(0)=1$.

$$\text{\boxed{\text{$Z_3'(A)=\Bigl[\frac{A(x)}{x},\frac{A(y)}{y},
    1-a, \frac{y-x}{y},l(y)\Bigr]$}}} $$
It's admissible because  $l(0)=1$ and $l((a-1)/a)=1$.

$$\text{\boxed{\text{$T_1\Big(\frac{A}{f}\Big)=\Bigl[\frac{A(x)}{x},\frac{A(y)}{y},
    \frac{(1-a)(x-1)}{x}, 1-\frac{x}{y}, l_1(y)\Bigr]$}}} $$
With the non-degeneracy assumption we have
\begin{alignat*}{3}
 \pa_1^0 T_1\Big(\frac{A}{f}\Big)&\subset \{t_3=1\},\quad&
 \pa_1^\infty T_1\Big(\frac{A}{f}\Big)&\subset \{t_4=1\},\quad
  \pa_2^\infty T_1\Big(\frac{A}{f}\Big)\subset \{t_5=1\},  \\
 \pa_3^\infty T_1\Big(\frac{A}{f}\Big)&\subset \{t_4=1\},\quad&
 \pa_4^\infty T_1\Big(\frac{A}{f}\Big)&\subset \{t_5=1\},   \quad
 \pa_5^\infty T_1\Big(\frac{A}{f}\Big)\subset\{t_4=1\},
\end{alignat*}
and
\begin{align*}
 \pa_2^0 T_1\Big(\frac{A}{f}\Big)&=\Bigl[\frac{A(x)}{x},
    \frac{(1-a)(x-1)}{x},\frac{aA(x)}{1-a},\frac{ac-c+1}{ac}\Bigr],\\
 \pa_3^0  T_1\Big(\frac{A}{f}\Big)&=\Bigl[\frac{1-a}{a},\frac{A(y)}{y},
       \frac{y-1}{y}, l_1(y) \Bigr],\\
 \pa_4^0  T_1\Big(\frac{A}{f}\Big)&=\Bigl[\frac{A(y)}{y},\frac{A(y)}{y},
    \frac{(1-a)(y-1)}{y}, l_1(y) \Bigr],\\
 \pa_5^0  T_1\Big(\frac{A}{f}\Big)&=\Bigl[\frac{A(x)}{x},\frac{A(c)}{c},
    \frac{(1-a)(x-1)}{x}, 1-\frac{x}{c} \Bigr].
\end{align*}
All these cycles are clearly admissible.

$$\text{\boxed{\text{$T_2\Big(\frac{A}{f}\Big)=\Bigl[\frac{A(x)}{x},\frac{A(y)}{y},
    \frac{(1-a)(x-1)}{x}, 1-\frac{x}{y}, l_2(y)\Bigr]$}}} $$
The above proof mostly is still valid because $l_2(0)=1$ except
that
\begin{align*}
 \pa_5^0  T_2\Big(\frac{A}{f}\Big)&=\Bigl[\frac{A(x)}{x},\frac{A(y_2)}{y_2},
    \frac{(1-a)(x-1)}{x}, 1-\frac{x}{y_2} \Bigr],\\
\pa_5^\infty T_2\Big(\frac{A}{f}\Big)&=\Bigl[\frac{A(x)}{x},-\mu,
    \frac{(1-a)(x-1)}{x}, B(x) \Bigr]
\end{align*}
which are both admissible since $B(0)=1$.

 From \eqref{h9} we have
\begin{equation*}
[h,h, \gd v,q_4,l_1] =\Bigl[x,y,\frac{abx+1}{abA(x)},
 \frac{y-x}{A(y)}, l_1(y)\Bigr].
\end{equation*}
It is clear that
$$\Bigl[x,y,\frac{abx+1}{abA(x)}, \frac{y}{A(y)}, l_1(y)\Bigr]$$
is admissible by $l_1(0)=1$. So we have
\begin{equation*}
[h,h, \gd v,q_4,l_1] =\Bigl[x,y,\frac{abx+1}{abA(x)},
 \frac{y-x}{y}, l_1(y)\Bigr]
\end{equation*}
which is also admissible.

\medskip
\noindent {\bf Step} (10). Final computation of $\{k(c)\}$.

\medskip

We have shown in the above everything in this step is admissible.
This completes the admissibility check of our main paper \cite{G}.

\end{document}